\documentclass{amsart}
\usepackage[utf8]{inputenc}
\usepackage[T1]{fontenc}
\usepackage{amssymb,amsfonts,amsthm,xcolor,microtype,siunitx,tikz,enumitem,calc,array,booktabs}
\usepackage[only,obslash]{stmaryrd}
\usepackage[hidelinks]{hyperref}
\setlist[itemize]{itemsep=1ex}
\numberwithin{equation}{section}
\newtheorem{prop}{Proposition}
\newtheorem{lemma}[prop]{Lemma}
\newtheorem{propdef}[prop]{Proposition and Definition}
\newtheorem{conjecture}[prop]{Conjecture}
\theoremstyle{definition}
\newtheorem{remark}[prop]{Remark}
\newtheorem{classlist}{Class List}

\newcommand{\E}[1]{\href{https://teorth.github.io/equational_theories/implications/?#1}{\textnormal{E}#1}}
\newcommand{\mkeqN}[1]{\overset{\!\!\E{#1}\!\!}=}

\newcommand{\opA}{*}
\newcommand{\opB}{/\!/}
\newcommand{\opC}{\backslash\!\backslash}
\newcommand{\copA}{\makebox[\widthof{$\opB$}]{\vphantom{$/$}$\opA$}}
\newcommand{\Qaxioms}{\mathfrak{Q}}
\newcommand{\ZZ}{\mathbb{Z}}
\newcommand{\FF}{\mathbb{F}}
\newcommand{\RT}{\mathrm{RT}}
\newcommand{\coloneqq}{:=}

\DeclareMathOperator{\Aut}{Aut}
\newcommand{\structure}[1]{\textbf{\boldmath #1}}
\newcommand{\file}[1]{\texttt{\detokenize{#1}}}

\begin{document}
\title{Implication semilattice \break of 990 quadratic quasigroup equational laws}
\author{Bruno Le Floch}
\date{November 2026}

\begin{abstract}
In his quest to exhibit a non-distributive lattice, Ernst Schr\"oder considered 135 years ago a list of \num{990} equational laws on quasigroups, a non-associative generalization of groups.  The laws of interest are analogues of associativity, such as $(x \opB y) \opA z = (y \opB x) \opC z$, built from the multiplication~$\opA$ and right/left conjugate-division operations, with the same three distinct variables on both sides.
I determine all \num{114} equivalence classes of their conjunctions, and all implications between them.  This includes as a small corner the non-distributive lattice identified by Schr\"oder.  Apart from \num{31} classes that are related to inverse properties (and weaker than total symmetry), and \num{3} that describe groups and their parastrophes, the remaining \num{80} are subject to the Four Quasigroups Theorem and quasigroups satisfying them are isotopes of Abelian groups.  More than half of these describe (affine) modules over interesting rings.
\end{abstract}

% arXiv comment: v1 9 pages. v2 37 pages with a description of each structure
% arXiv title: Implication semilattice of 990 quadratic quasigroup equational laws
% arXiv author: Bruno Le Floch
% arXiv abstract: In his quest to exhibit a non-distributive lattice, Ernst Schr\"oder considered 135 years ago a list of 990 equational laws on quasigroups, a non-associative generalization of groups.  The laws of interest are analogues of associativity, such as $(x // y) * z = (y // x) \backslash\backslash z$, built from the multiplication $*$ and right/left conjugate-division operations, with the same three distinct variables on both sides. I determine all 114 equivalence classes of their conjunctions, and all implications between them.  This includes as a small corner the non-distributive lattice identified by Schr\"oder.  Apart from 31 classes that are related to inverse properties (and weaker than total symmetry), and 3 that describe groups and their parastrophes, the remaining 80 are subject to the Four Quasigroups Theorem and quasigroups satisfying them are isotopes of Abelian groups. More than half of these describe (affine) modules over interesting rings.
%
% files to upload:
% Schroeder-990.tex
% Schroeder-990.bbl
% ATP_utils_quasigroup.py
% Schroeder_equations.py
% Schroeder_implications.py
% Schroeder_conjunctions.py
% Schroeder_find_class.py
% quasigroup-equations.txt
% quasigroup-representatives.txt
% quasigroup-classes.txt
% quasigroup-implications.txt
% quasigroup_long_classes.py
% quasigroup_long_classes.json
% quasigroup_eq_to_long_class.json
% Sch-102.p

\maketitle

\tableofcontents

\section{Motivation and main results}

Riding on the momentum of the Equational Theories Project~\cite{ETP}, this report pays tribute to seemingly the first attempt\footnote{As a non-expert, I welcome corrections on this point, as well as terminology and notation.  Throughout this work I will use modern terminology, with no attempt at a historical account of where the words originated.  This is not motivated by historical revisionism but rather by lack of expertise in the history of mathematics.} to wrestle with a large set of equational laws, by completing a classification work begun by Ernst Schr\"oder in the 1890's, explained in his book \cite[Anhang 4--5, pp. 617--646]{Schroder-1890}.  Schr\"oder developed the theory of quasigroups, which are sets equipped with a binary operation~$\opA$ (not necessarily associative) whose right and left multiplications are bijective.

The theory of quasigroups is vast (see e.g., the textbook~\cite{pflugfelder1990quasigroups}), and it is impossible to properly account for the whole history of the field in just a few lines.  Throughout the decades, researchers have considered a multitude of properties that quasigroups can satisfy, such as inverse properties (of relevance in \S\ref{sec:twovars}) or Bol--Moufang-type laws, structures associated to a given quasigroups (such as their translation groups, center, nucleus, etc.\@), and relations between quasigroups (such as isotopy).

The rest of this section is devoted to quickly presenting the problem and results.
Then \S\ref{sec:tools} describes the relation to the Equational Theories Project, the tools I used, and how to use the resulting data files.
Finally, \S\ref{sec:twovars} and subsequent sections describe all $114$ quasigroup varieties of interest, in terms of known quasigroup properties, which are reviewed as needed.
Appendices collect a few asides: an example quasigroup~\S\ref{app:counter};
asymptotic features of the number of quasigroups for certain classes~\S\ref{sec:pressure}, and figures with several sublattices~\S\ref{sec:lattice-figs}.

\subsection{Equational laws and quasigroup varieties}
\label{sec:eqlaws}

A quasigroup is a set~$Q$ equipped with a binary operation~$\opA$ whose right and left multiplications are bijective.  Inverting these maps defines right and left division operations $/$ and~$\backslash$ characterized by $(x/y)\opA y=x$ and $y \opA (y\backslash x)=x$.  It is more symmetrical to work with their conjugates $x\opB y = y/x$ and $x\opC y = y \backslash x$.  From this point of view, a quasigroup is a set $Q$ equipped with three binary operations $\opA,\opB,\opC\colon Q\times Q\to Q$ such that for all $x,y\in Q$ one has
\begin{equation}\label{quasigroup-def-3}
  \Qaxioms : \left\{ \quad
  \begin{aligned}
    x \opB (y \opA x) & = y , \qquad & (x \opB y) \opA x & = y , \\
    x \opA (y \opC x) & = y , \qquad & (x \opA y) \opC x & = y , \\
    x \opC (y \opB x) & = y , \qquad & (x \opC y) \opB x & = y .
  \end{aligned}
  \right.
\end{equation}
These identities can also be summarized as the equivalence
\begin{equation}\label{xyz-iff}
  \forall x,y,z\in Q, \qquad \bigl(x \opA y = z \iff y \opB z = x \iff z \opC x = y\bigr) .
\end{equation}
The resulting definition of quasigroups~\eqref{quasigroup-def-3} features a $\ZZ/3\ZZ$ symmetry under which one cyclically permutes the operations $(\opA,\opB,\opC)$.  This symmetry extends to an $S_3$~permutation symmetry called \emph{parastrophy}, which maps the quasigroup $(Q,\opA)$ to one of the six quasigroups $(Q,\opA)$, $(Q,\tilde\opA)$, $(Q,\opB)$, $(Q,/)$, $(Q,\opC)$, $(Q,\backslash)$, where $\tilde\opA$ is defined by $x\,\tilde\opA\,y=y\opA x$.

Following Schr\"oder, I focus on \num{990} equational laws of the form
\begin{equation}\label{Schroeder-form}
  \forall x, \forall y, \forall z, \quad x \odot_1 y \odot_2 z = \alpha \odot_3 \beta \odot_4 \gamma
\end{equation}
where $\odot_i$ are replaced (independently) by any of the three quasigroup operations, $\alpha,\beta,\gamma$ are the variables $x,y,z$ in some order, and the operations are associated either as $\_\odot(\_\odot\_)$ or $(\_\odot\_)\odot\_$ independently on both sides.  One modern terminology is to call these balanced quadratic equations in three variables on quasigroups~\cite{sokhatsky2024quadratic}; another is to dub them linear equations~\cite{niemenmaa1992general}.

Any set~$A$ of equational laws defines a variety $V(A)$ consisting of quasigroups\footnote{I shall be sloppy (from a universal algebra point of view) in identifying quasigroups $(Q,\opA,\opB,\opC)$ with three operations, which are type-$(2,2,2)$ algebras, to quasigroups $(Q,\opA)$ with a single operation subject to the condition that left and right multiplications be bijective.} satisfying every law in~$A$.
The quasigroups automatically satisfy any logical consequence of the laws in~$A$.
It is thus natural to focus on subsets of the \num{990} equations that are closed under logical implications.
One can define a lattice structure by introducing two binary operations on pairs of such subsets $A$ and~$B$.
\begin{itemize}

\item The logical sum $A+B$, which contains the set-theoretic union $A\cup B$, is the set of all consequences of $A\land B$ among the \num{990} laws.  This operation is natural at the level of quasigroup varieties since $V(A+B) = V(A) \cap V(B)$.

\item The product $AB$ is the set-theoretic intersection $A\cap B$.  In terms of varieties one has $V(AB) \supset V(A) \cup V(B)$, with no equality in general as $AB$ may be a very meager set of laws.

\end{itemize}
Another operation (which can be reconstructed from either of these) is implication: for subsets closed under logical implication, $A\implies B$ if $B\subset A$; equivalently $A+B=A$ and equivalently $AB=B$.
Since the product operation is less natural for quasigroup varieties, I find it more proper in the title to use the terminology \emph{semilattice} (equipped with the operation~$+$) rather than \emph{lattice}.

\subsection{Schr\"oder's non-distributive lattice}

Schr\"oder studied a few such subsets: the set~$A_1$ of $16$~consequences of associativity (Sch-253\,$\Longleftrightarrow$\,Sch-65)\footnote{Numbers of the form Sch-$NN$ refer to an arbitrary numbering of the Schr\"oder laws, described in \S\ref{subsec:numbering}.  The associativity equation $x \opA (y \opA z) = (x \opA y) \opA z$ (Sch-253) is equivalent to $x \opA (y \opB z) = y \opB (x \opA z)$ (Sch-65), which is how it is listed later on due to the ordering.}, the set $C_1$ of $30$~consequences of commutativity of~$\opA$ (Sch-1), the set $O_1=A_1+C_1$ of $150$~consequences of associativity and commutativity (Sch-3), the set $C_{00}$ of $18$~consequences of $x \opA (y \opA z) = x \opB (y \opB z)$ (Sch-23), and he showed that
\begin{equation}
  \begin{aligned}
    A_1C_{00} + C_1C_{00} & = \emptyset + \emptyset = \emptyset , \\
    E_1 \coloneqq (A_1 + C_1)C_{00} & = \bigl\{\textnormal{Sch-127}, \textnormal{Sch-899}\bigr\} \neq \emptyset ,
  \end{aligned}
\end{equation}
(where Sch-127 is $x \opA (y \opC z) = x \opC (y \opB z)$ and Sch-899 is equivalent to it).  These sets differ so the distributive law fails in the lattice.
Schr\"oder also points out several other interesting subsets, whose notation I include in \S\ref{sec:twovars} and subsequent sections.  For instance, $Q_0 = C_1+C_{00}$ and $Q_1 = O_1 Q_0 = C_1+O_1C_{00}$ and $K_1=A_1Q_0$.

In retrospect, the resulting lattice could be described simply by its Hasse diagram (see Figure~\ref{fig:five-lattice}) with no reference to quasigroups, but the viewpoint of prescribing a lattice as some given set equipped with suitable operations given by explicit tables was only understood later on by Dedekind.
Even though the original lattice motivation is gone, the question of exploring varieties of quasigroups remains, and I completed Schr\"oder's unfinished quest of determining the whole semi-lattice for the \num{990} equational laws.

\begin{figure}
  \centering
  \begin{tikzpicture}
    \node[draw](all) at (0,0) {$\emptyset$: All quasigroups};
    \node[draw](sch336) at (-4.8,-1) {$K_1$: EIP (Sch-336)};
    \node[draw](mip) at (4.4,-1) {$E_1$: MCIP (Sch-127)};
    \node[draw](com) at (0,-1) {$C_1$: Commutation (Sch-1)};
    \node[draw](gp) at (-5,-2) {$A_1$: Groups (Sch-65)};
    \node[draw](c00) at (6,-2) {$C_{00}$: ECIP (Sch-23)};
    \node[draw](tts) at (0,-3) {$Q_1$: $\ZZ/2\ZZ$-twisted total symmetry (Sch-88)};
    \node[draw](ab) at (-4.5,-4) {$O_1$: Abelian groups (Sch-3)};
    \node[draw](tot) at (4.5,-4) {$Q_0$: Total symmetry (Sch-24)};
    \node[draw](bool) at (0,-5) {$A_1+C_{00}$: Boolean groups (Sch-7)};
    \draw[-{stealth}](mip)--(all);
    \draw[-{stealth}](tts)--(mip);
    \draw[-{stealth}](sch336)--(all);
    \draw[-{stealth}](gp)--(sch336);
    \draw[-{stealth}](tts)--(sch336);
    \draw[-{stealth}](com)--(all);
    \draw[-{stealth}](c00)--(mip);
    \draw[-{stealth}](tts)--(com);
    \draw[-{stealth}](tot)--(c00);
    \draw[-{stealth}](ab)--(gp);
    \draw[-{stealth}](tot)--(tts);
    \draw[-{stealth}](ab)--(tts);
    \draw[-{stealth}](bool)--(tot);
    \draw[-{stealth}](bool)--(ab);
  \end{tikzpicture}
  \caption{\label{fig:five-lattice}Sublattice (of the full lattice of logically closed subsets of \num{990} laws) generated by Schr\"oder's counterexample to distributivity.  Arrows describe the Hasse diagram of the implication poset; the logical sum and products are meet and join in this poset.  MCIP is the middle crossed inverse property, and EIP (ECIP) denote the equal (crossed) inverse property, see \S\ref{sec:twovars}.}
  % python Schroeder_find_class.py hasse 1 65 23
\end{figure}

\subsection{The Four Quasigroups Theorem}
\label{subsec:four-quasigroups}

Up to renaming variables, all \num{990} equations can be brought to one of the two forms (whose name will become clear momentarily)
\begin{equation}
  \begin{aligned}
    x \odot_1 (y \odot_2 z) & = x \odot_3 (y \odot_4 z) \qquad & & \text{(reducible),} \\
    x \odot_1 (y \odot_2 z) & = (x \odot_3 y) \odot_4 z \qquad & & \text{(generalized associativity),}
  \end{aligned}
\end{equation}
where $\odot_i$ are among the \emph{six} quasigroup operations, namely $\opA$, $\opB$, $\opC$ and the operations with their arguments swapped.
The \num{324} equations that take the first form, and their conjunctions, are studied in \S\ref{sec:twovars}.  As explained there, these equations are equivalent to two-variable equations, hence are reducible in the sense of~\cite{sokhatsky2024quadratic}.  They are intermediate between some well-known (crossed) inverse properties and total symmetry.

The \num{666} equations of the second kind take the shape of associativity with four possibly distinct operations.  They fall under the purview of the powerful Four Quasigroups Theorem\footnote{According to~\cite{dudek2005remarks}, this theorem was announced in 1958 in~\cite{belousov1958associative}, proven in~\cite{belousov1961globally}, with a modified proof presented in the book~\cite{aczel1963generalized}.  See also~\cite{krapez1980generalized}.  I have also found \cite[p.19]{shcherbacov2009structure} helpful.  An early version of this theorem with three of the operations being equal appeared already in 1929 in~\cite{Suschkewitsch:1929}.  We give a proof in~\S\ref{subsec:prove-four-quasigroups}.  Amusingly, an analogous result holds for other equations than associativity, such as the ``quadratical property'' $(x y) x = (zx) (yz)$~\cite{kovali2005functional}.}, which states that the quasigroup~$Q$ carries a group structure $(Q,\times)$ and
\begin{equation}\label{four-quasigroups-result}
  \begin{aligned}
    x \odot_1 y & = \alpha(x) \times \mu(y) , \qquad &
    x \odot_2 y & = \mu^{-1}(\beta(x) \times \gamma(y)) , \\
    x \odot_3 y & = \nu^{-1}(\alpha(x) \times \beta(y)) , \qquad &
    x \odot_4 y & = \nu(x) \times \gamma(y) ,
  \end{aligned}
\end{equation}
where $\alpha,\beta,\gamma,\mu,\nu$ are bijections of~$Q$.  In the present context, the operations $\odot_i$ are not independent since they are among the six operations of a given quasigroup.  Depending on which of the \num{666} equations is considered, this constrains the bijections.  In all cases of interest to us (except for usual associativity itself and its parastrophes), Belousov's criterion ---which we shall state in \S\ref{subsec:Belousov-criterion}--- leads to $(Q,\times)$ being an Abelian group, cf.\ Proposition~\ref{prop:trichotomy}.

\emph{There remains two puzzles.}
\begin{itemize}
\item For each non-reducible equation considered here that is not equivalent to associativity or its parastrophes (described in~\S\ref{sec:groups}), there is a way to apply Belousov's criterion to conclude that the underlying group is Abelian, as proven in Proposition~\ref{prop:trichotomy}, but there is no systematic understanding of why this is always the case.

\item Many classes imply linearity and there is no clear description of why, analogous to the Belousov criterion.
\end{itemize}

\subsection{Main result}
\label{subsec:main-result}

The main outcome of this work consists of a classification of all subsets of equations closed under implication.  These \num{114} subsets are given in \file{quasigroup_long_classes.json} (see details in \S\ref{subsec:supplements}).  Given two such subsets $A$ and~$B$, the logical operations of interest are determined as follows:
$A+B$ is the smallest such subset containing both $A$ and~$B$;
$AB$ is the intersection of $A$ and $B$, which is automatically logically closed;
$A\implies B$ if and only if $B\subset A$.
Given \emph{any subset} of the \num{990} equations, the closure under logical implications is the smallest subset listed in the data files that contains all of these equations.
Some interesting sublattices are depicted in figures throughout the paper.

The \num{114} equivalence classes of conjunctions (or equivalently, the quasigroup varieties) are usefully organized according to how constraining they are.
They are listed in \S\ref{sec:twovars}--\ref{sec:isotopy} together with a description in conventional mathematical terms.
\begin{itemize}
\item[\S\ref{sec:twovars}] describes the \num{31} (classes of) conjunctions of reducible equations, understood as twists of commutativity, or semi-symmetry, or total symmetry, by a bijection that is typically an automorphism of small order.
\item[\S\ref{sec:groups}] describes \num{3} classes that state associativity for each of the three operations.
\item[\S\ref{sec:isotopy}] explains that quasigroups in any of the remaining classes carry an Abelian group structure (modulo a choice of origin).
\item[\S\ref{subsec:quasilinear}] lists \num{13} classes describing Abelian groups equipped with bijections that are not necessarily linear maps.
\item[\S\ref{subsec:noncomm}] lists \num{19} classes of non-commutative linear models (over Abelian groups).
\item[\S\ref{subsec:medial}] lists \num{48} classes of commutative linear models, also called medial quasigroups, which describe modules or affine modules over specific rings.
\end{itemize}
Some of the varieties were already identified as being interesting in the context of the Equational Theories Project (ETP), such as affine modules over $\ZZ[\sqrt{-1}]$.
Many classes are parastrophic to earlier ones, in the sense that they are transformed to them upon permuting the six quasigroup operations; the commentary then mostly points to the other class using the abbreviation ``para''.
For some of the classes, I noted the equivalence to a law written only in terms of one binary operation, in which case I give the ETP number as E\textit{nnn} with a link to the \href{https://teorth.github.io/equational_theories/implications}{interactive Equation Explorer} that provides more commentary in some cases.

It would be interesting to extend the semi-lattice with other collections of equations, such as other quadratic equations in~$\opA$ only~\cite{krapevz2007quadratic,akhtar2016generalized,akhtar2017symmetric,akhtar2022linear} or equations of Bol--Moufang type \cite{phillips2005varieties,phillips2005loopvarieties,akhtar2012varieties}.

\subsection{Data availability statement}

Data files and scripts produced in the course of this research are available as ancillary files at \url{https://arxiv.org/abs/2603.29909}  Their use is described in \S\ref{subsec:supplements}.

\subsubsection{AI use statement}

When this research was performed (early 2026), freely available AI were not so useful.  I made some unsuccessful attempts to use an LLM to help figure out standard terminology around crossed inverse properties, but reading articles was more efficient, as was reading software documentation (Python, Mace4, Prover9) for writing the relevant scripts.

\section{Tools used}
\label{sec:tools}

\subsection{Numbering and counting Schr\"oder's laws}
\label{subsec:numbering}

Since the original work did not appear to specify an ordering of the equations, I decided arbitrarily on a numbering, prefixed with ``Sch-'' in reference to Schr\"oder.\footnote{I plan to keep this numbering fixed in any future version of this work despite its shortcomings.  To allow for an easier migration to a more general numbering of all quasigroup equational laws, I have included the prefix ``Sch-'' throughout the discussion.}  Namely, as can be seen in the data file \file{quasigroup-equations.txt}, equations are ordered
\begin{itemize}
\item first by how the two binary operations on the left and right are associated, namely $\_\odot(\_\odot\_) = \_\odot(\_\odot\_)$ before $\_\odot(\_\odot\_) = (\_\odot\_)\odot\_$ before $(\_\odot\_)\odot\_ = (\_\odot\_)\odot\_$;
\item then by the choice of the four operations~$\odot$ in lexicographical order, namely $({\opA},{\opA},{\opA},{\opA})$, then $({\opA},{\opA},{\opA},{\opB})$, then $({\opA},{\opA},{\opA},{\opC})$, then $({\opA},{\opA},{\opB},{\opA})$, etc.;
\item then by the choice of variables.
\end{itemize}
Following Schr\"oder, and following what we did in the Equational Theories Project~\cite{ETP}, all tautological equations are omitted, as are any equation that can be obtained from earlier ones by relabelling variables or swapping left and right hand side.  This explains the absence of any equation associated as $(\_\odot\_)\odot\_ = \_\odot(\_\odot\_)$.

The left-hand side always uses the variables $x,y,z$ in this order, while they can appear in any of the six orders on the right-hand side, unless both sides involve the same operations associated in the same way, in which case one may not take $(\alpha,\beta,\gamma)$ to be $(x,y,z)$ in the same order (tautology) nor $(z,x,y)$ because this is equivalent to the choice $(y,z,x)$.

This gives $9^2\times 6$ laws associated as $\_\odot(\_\odot\_)=(\_\odot\_)\odot\_$,
$\binom{9}{2}\times 6$ of the form $\_\odot(\_\odot\_)=\_\odot(\_\odot\_)$ with different pairs of operations on the two sides and $9\times 4$ of that form with the same pair of operations on both sides, and then the same number $\binom{9}{2}\times 6+9\times 4$ with the dual association $(\_\odot\_)\odot\_=(\_\odot\_)\odot\_$, for a total number of equations
\begin{equation}
  9^2\times 6 + 2 \times \binom{9}{2}\times 6 = 990 .
\end{equation}

\subsection{The Equational Theories Project}

The Equational Theories Project~\cite{ETP} considered magmas, which are sets equipped with one binary operation, and the $4694$ simplest equational laws, that is, all laws that involve at most~$4$ operations.
This collaborative research project determined all implications between these laws.\footnote{In the context of finite magmas, one implication remains undetermined: whether in all finite magma such that $\forall x,\ \forall y,\  x = y \diamond (x \diamond ((y \diamond x) \diamond y))$ one necessarily has $\forall x, \ x = ((x\diamond x)\diamond x)\diamond x$.}

Automated theorem provers were particularly effective tools to obtain positive implications: in all cases where an equation $E$ implies another equation~$E'$, a short proof could be found automatically.
For counterexamples, the situation was somewhat trickier: small countermodels could be obtained directly by automated tools, but larger finite countermodels could only be found by imposing well-chosen assumptions to restrict the search.  Even worse, some implications turned out to hold for all finite magmas but not for infinite magmas.  To disprove such implications it was necessary to build infinite countermodels.
Thus, overall, the Equational Theories Project relied on a collection of techniques:
\begin{itemize}
\item Automated theorem provers and SAT solvers such as Vampire, Prover9, Mace4.
\item A Lean formalization of every implication and every countermodel, to ensure nothing was missed.  This was absolutely crucial due to the large number of collaborators and the tens of millions of implications to be considered.
\item Mathematical developments, to come up with new constructions of magmas based on magma cohomology, greedy constructions, and rewriting theory.
\end{itemize}
(As well as tools to organize the work, visualize the status of the project, explore the implication graph, and much more.)

\subsection{Methods}

The main goal of the present project is to determine all logically-closed subsets~$A$ of equations among Schr\"oder's \num{990} equations.  This relied on two automated tools developed by McCune: the automated theorem prover (ATP) Prover9, which given some equations and a goal seeks a proof thereof, and the model builder Mace4, which given the same setup seeks a finite countermodel, namely (in our case) a quasigroup that does not obey the target equation.
This was combined with simple Python scripts.

For any given implication of interest between a set~$A$ of equations and a target~$E$, a useful strategy was as follows.  (Here we denote by $\Qaxioms$ the quasigroup axioms~\eqref{quasigroup-def-3} above.)
Seek countermodels to $\Qaxioms\land A\implies E$ using Mace4 with a very small bound in size.
Check if the implication can be proven using Prover9, with a one-second timeout.
Finally make a longer attempt at counterexamples using Mace4.\footnote{Some fiddling with Mace4 parameters such as \texttt{selection\_measure} was needed to get a fast enough answer.}
This always succeeds for the implications of interest, hence serves as an \emph{oracle} for the truth of implications.
This is a crucial simplification compared to the Equational Theories Project: all of the non-implications admitted a \emph{finite counterexample}, which could be found automatically.

The strategy to find all inequivalent conjunctions of equations is thus as follows.
\begin{itemize}
\item For each pair of equation $E,E'$ (among the \num{990}), determine whether $\Qaxioms\land E\implies E'$.
\item Deduce the \num{47} lowest-numbered representatives of equivalence classes of equations.  The consequences of a given equation form a logically closed subset.
\item For each logically closed subset~$A$ found so far, in turn, and each of the $47$ representatives~$E$ that is not in~$A$, consider $A\cup\{E\}$.  Seek all implications of the form $\Qaxioms\land A\land E\implies E'$ for $E'$ among the \num{47} representatives.  Many of these searches are skipped thanks to transitivity of implication.  If the resulting set of consequences is new, add it to the list of logically closed subsets to consider.
\end{itemize}
One should remember as well that the empty set is logically closed!

\subsection{Code and data files}
\label{subsec:supplements}

Supplementary files are accessible at the address \url{https://arxiv.org/abs/2603.29909}.
Python scripts used in this work are as follows.  Some require the scripts \texttt{mace4} and \texttt{prover9} to be in the same directory.
\begin{itemize}
\item \file{ATP_utils_quasigroup.py} provides an interface to the Mace4 and Prover9 ATP tools, which must be installed separately.
\item \file{Schroeder_equations.py} generates the list of \num{990} equations in the three operations $\opA,\opB,\opC$.
\item \file{Schroeder_implications.py} determines all implications (and equivalences) between these equations, by running Mace4 and Prover9.
\item \file{Schroeder_conjunctions.py} assembles equations between the $47$ representatives and finds all inequivalent conjunctions.
\item \file{Schroeder_find_class.py} converts from equations (possibly not with canonical variable names) to equation numbers, finds the class containing a list of equations, finds all classes weaker/stronger than a class, and generates Hasse diagrams; for instance the following code served as a starting point for \autoref{fig:five-lattice}:
\begin{verbatim}
    python Schroeder_find_class.py hasse 1 65 23
\end{verbatim}
\end{itemize}

The scripts generate the following data files, also included.
\begin{itemize}
\item \file{quasigroup-equations.txt} gives the list of \num{990} equations.
\item \file{quasigroup-representatives.txt} gives the \num{47} representatives of equivalence classes.
\item \file{quasigroup-classes.txt} gives each representative together with the corresponding list of equation numbers in the equivalence class.
\item \file{quasigroup-implications.txt} gives all implications between representatives: all other possible implications are disproven by some finite quasigroup.
\item \file{quasigroup_long_classes.py} gives the list of long classes, namely of logically-closed subsets of equations, ordered by cardinality, as well as the mapping between representatives and long classes.  The same data is given in \file{quasigroup_long_classes.json} and \file{quasigroup_eq_to_long_class.json} in the JSON format.
\end{itemize}

As an example of Prover9 file, \file{Sch-102.p} is also provided; it shows Proposition~\ref{prop:102-boolean}.

\section{Consequences of total symmetry (31 classes)}
\label{sec:twovars}

\subsection{A first type of equation}

Among the \num{990} equations, \num{324} can be written as
\begin{equation}\label{xyz-xyz}
  x \odot_1 (y \odot_2 z) = x \odot_3 (y \odot_4 z)
\end{equation}
where $\odot_i$ are among the six quasigroup operations, namely $\opA,\opB,\opC$ and their conjugates $\tilde\opA,/,\backslash$ (which have swapped arguments).
They are equivalently characterized as the consequences of total symmetry (Sch-24), which states that all six quasigroup operations are equal.
Up to equivalence, these equations reduce to only \num{21} and their conjunctions fall into \num{31} equivalence classes.

Section \S\ref{subsec:reducible-eqs} explains how to reduce equations of the form~\eqref{xyz-xyz} to single-variable equations and equalities between translation sets.
The latter amount to known properties of quasigroups such as inverse or crossed inverse properties.
A useful organizing principle is how many of the six quasigroup operations are identified (up to twists) by the equations.
This principle is presented in \S\ref{subsec:fine-spec}, and put on a sharper footing in terms of the asymptotic number of quasigroups of a given cardinality.
This leads to presenting the \num{31} classes in \S\ref{subsec:twist-comm}--\S\ref{subsec:twist-ts}:
\begin{itemize}
\item \num{1} class has no equation and corresponds to the variety of all quasigroups, see \autoref{class:empty} in \S\ref{subsec:twist-comm};
\item \num{6} classes relate the six operations pairwise: commutativity and (left/middle/\allowbreak right) inverse property (Definition~\ref{propdef:lmr-inverse-prop}), see \autoref{class:comm} in \S\ref{subsec:twist-comm};
\item \num{10} classes relate the three operations $\opA,\opB,\opC$ together (and relate $\tilde\opA,/,\backslash$): semi-symmetry and its consequences, such as crossed inverse properties (Definition~\ref{propdef:lmr-crossed-inverse-prop}), see \autoref{class:ss} in \S\ref{subsec:twist-ss};
\item \num{14} classes relate all six operations: total symmetry and twists thereof, see \autoref{class:ts} in \S\ref{subsec:twist-ts}.
\end{itemize}
Within these sections I also suggest some names for the quasigroup varieties that do not appear to have been named in the literature on quasigroups, such as (left/middle/right) \emph{involutive crossed inverse properties}, the \emph{equal inverse property}, etc.\@
The implication semi-lattice between these classes is depicted in \autoref{fig:two-vars-lattice}.

\subsection{Reducible equations}
\label{subsec:reducible-eqs}

Equations of the form~\eqref{xyz-xyz} are \emph{reducible equations} in the sense that they are equivalent to equations with fewer variables.  They are conveniently expressed as equalities between sets of translations, as stated in Proposition~\ref{prop:reduce}.  For some quasigroup operation~$\odot_i$, the right-translation by some $w\in Q$ is the bijection $x\mapsto x\odot_i w$, and the set of these right-translations is
\begin{equation}
  \RT(\odot_i) = \bigl\{ x \mapsto x \odot_i w , \ w\in Q\bigr\} .
\end{equation}
We also need the notion of left division~$\obslash_i$, defined by $x \obslash_i (x \odot_i y) = y$.

\begin{prop}\label{prop:reduce}
  The equational law $x \odot_1 (y \odot_2 z) = x \odot_3 (y \odot_4 z)$ is equivalent to a conjunction of two-variable equations.  It is also equivalent to $\RT(\odot_1) = \RT(\odot_3)$ and $\RT(\odot_2) = \RT(\odot_4)$ and the equational law $x \obslash_3 (x \odot_1 x) = x \odot_4 (x \obslash_2 x)$.

  In the special case where $\odot_1=\obslash_2$ and~$\odot_3=\obslash_4$, the equational law is simply equivalent to $\RT(\odot_1) = \RT(\odot_3)$.
\end{prop}

This result, proven momentarily, relates the \num{324} reducible equations to equalities of translation sets which were considered recently by Sokhatsky and Lutsenko~\cite{sokhatsky2020classification,sokhatsky2021classification}.  These authors found the \num{9} quasigroup varieties corresponding to the ``special case'' in Proposition~\ref{prop:reduce}.  These are
\begin{itemize}
\item the left/middle/right crossed inverse property (\num{3} varieties);
\item the left/middle/right inverse property (\num{3} varieties);
\item commutativity\footnote{The authors of~\cite{sokhatsky2020classification,sokhatsky2021classification} defined ``mirror quasigroups'' as obeying the identity $\varphi(x)\opA y = y\opA x$ for some function~$\varphi$.  By evaluating this identity with $y=x$ one gets $\varphi(x)\opA x=x\opA x$, and right-dividing by~$x$ gives $\varphi(x)=x$, so that a mirror quasigroup is nothing but a commutative quasigroup.} and its parastrophes (\num{3} varieties).
\end{itemize}
Conjunctions of these equations define \num{22} classes denoted ``no eq'', MCIP, CIP, MIP, IP, all IP, ``commutes'', 4-TS, 2-TS, and TS below (see \autoref{fig:two-vars-lattice}).  The remaining \num{9} classes EIP, MICIP, ICIP, ECIP, SS involve not only the equality of several translation sets, but also single-variable identities of the form $x \obslash_3 (x \odot_1 x) = x \odot_4 (x \obslash_2 x)$.

\begin{proof}[Proof of Proposition~\ref{prop:reduce}]
By changing variables to $w = y \odot_2 z$ and acting with $x\obslash_3\_$ on both sides one gets $x \obslash_3 (x \odot_1 w) = y \odot_4 (y \obslash_2 w)$.  The resulting equation is equivalent (for instance) to a collection of equations in one and two variables,
\begin{equation}\label{three-two-vars}
  x \obslash_3 (x \odot_1 w)
  = w \obslash_3 (w \odot_1 w)
  = w \odot_4 (w \obslash_2 w)
  = y \odot_4 (y \obslash_2 w) .
\end{equation}
It is natural to call these \num{324} equations ``reducible'' as a result of this reduction in the number of variables.  Equivalently,
\begin{equation}\label{less-two-vars}
  x \obslash_3 (x \odot_1 w) = f(w) ,
  \qquad
  y \odot_4 (y \obslash_2 w) = f(w) ,
\end{equation}
for the \emph{same} function $f:Q\to Q$.  This function is necessarily bijective since it can be expressed (for any choice of $x\in Q$) as the composition of the left translations $w \mapsto x\odot_1 w$ and $w\mapsto x \obslash_3 w$.

The first identity in~\eqref{less-two-vars} can be rewritten as $x \odot_1 w = x \odot_3 f(w)$, namely as the equality of the right translations $x\mapsto x \odot_1 w$ and $x \mapsto x \odot_3 f(w)$.  Such a bijection~$f$ exists if and only if the set of right translations $\RT(\odot_1)$ coincides with~$\RT(\odot_3)$.
In the same way, the second identity in~\eqref{less-two-vars} is equivalent to $\RT(\odot_2)=\RT(\odot_4)$.  There remains to state that the bijections coincide.  This is precisely the identity $x \obslash_3 (x \odot_1 x) = x \odot_4 (x \obslash_2 x)$.

The special case is obvious since the two equalities of right translation sets become identical, and the extra one-variable identity also trivializes when $\odot_1,\odot_3$ are left divisions of $\odot_2,\odot_4$, respectively.
\end{proof}

\subsection{Comment on the asymptotic fine spectrum}
\label{subsec:fine-spec}

The \num{324} equations considered here are also characterized as consequences of the totally symmetric quasigroup law Sch-24 (listed below), namely they hold whenever $\opA$, $\opB$, $\opC$ are all equal and commutative.
As stated in (the proof of) Proposition~\ref{prop:reduce}, the equations and their \num{31} conjunctions amount to various equalities of the form $x \odot_1 y = x \odot_3 f(y)$ for a bijection~$f$.  Such an equality induces analogous relations between the left and right divisions $\oslash_i,\obslash_i$ of the two operations,
\begin{equation}\label{xobslash3y}
  x \oslash_3 f(y) = x \oslash_1 y , \qquad
  x \obslash_3 y = f(x \obslash_1 y) .
\end{equation}
All of these relations take the general form $x \odot y = f(g(x) \odot' h(y))$ for three bijections $f,g,h$, namely they state that $\odot,\odot'$ are \emph{isotopic}.
What is found in~\eqref{xobslash3y} is the standard result that $\odot_1,\odot_3$ are isotopic if and only if $\oslash_1,\oslash_3$ are isotopic, if and only if $\obslash_1,\obslash_3$ are isotopic.

Depending on which pairs of operations are isotopic, each of the \num{31} classes of interest fall into a few basic cases (besides the class with no equation at all).
\begin{itemize}
\item If $\opA$ is isotopic to $\tilde{\opA}$, then $\opB$ is isotopic to~$\backslash$, and also $\opC$~to~$/$.  Likewise if $\opA$ is isotopic to~$\backslash$ or to~$/$, in which we get isotopy of $\opB$ with its conjugate or $\opC$ with its conjugate, respectively.  Such classes describe variants of commutativity, twisted by a bijection, discussed in \S\ref{subsec:twist-comm}.

\item If $\opA$ is isotopic to $\opB$ then they are also isotopic to~$\opC$, and the three conjugate operations $\tilde\opA,/,\backslash$ are themselves isotopic to each other.  Such classes describe variants of semi-symmetry, twisted by a bijection, discussed in \S\ref{subsec:twist-ss}.

\item If both isotopies hold simultaneously, then all six operations are isotopic.  Such classes describe variants of total symmetry, twisted by a bijection, discussed in \S\ref{subsec:twist-ts}.
\end{itemize} 
I conjecture that these three cases are distinguished by the asymptotic number $N_n(A)$ of quasigroups of cardinality~$n$ in the given quasigroup variety~$V(A)$, as $n\to +\infty$.

\begin{conjecture}\label{conj:NnA}
  Each class~$A$ of conjunctions of reducible equations has
  \begin{equation}\label{NnA-2var}
    N_n(A) = n^{(1-p_A)n^2} e^{O(n^2)} = e^{(1-p_A)n^2\log n+O(n^2)}
  \end{equation}
  for some ``pressure'' $p_A \in \{0, 1/2, 2/3, 5/6\}$ that measures which proportion of the Cayley table of~$\opA$ is constrained by the equation.
  Specifically, $p_A=0$ for the empty class, $p_A = 1/2$ for twists of commutativity, $p_A = 2/3$ for twists of semi-symmetry, and $p_A = 5/6$ for twists of total symmetry.
\end{conjecture}

\begin{figure}
  \centerline{%
  \begin{tikzpicture}
    \filldraw[thick,fill opacity=.3, color=blue](-5.5,0) rectangle (-2,9);
    \filldraw[thick,fill opacity=.3, color=green](-1,-3)--(-1,5.8)--(3.7,5.8)--(3.7,1)--(1,0.1)--(1,-3)--cycle;
    \filldraw[thick,fill opacity=.3, color=orange](4,1.7) rectangle (6,8.9);
    \node[draw](Sch-) at (0,9.5) {\rotatebox{90}{no eq}};

    \node[draw](Sch-127) at (-3.6,8) {\rotatebox{90}{$\copA$ MCIP}};
    \node[draw](Sch-432) at (4.4,8) {\rotatebox{90}{$\copA$ MIP}};
    \node[draw](Sch-127-432) at (-.6,5) {\rotatebox{90}{$\copA$ 4-TS}};
    \node[draw](Sch-93) at (-3.6,5) {\rotatebox{90}{$\copA$ MICIP}};
    \node[draw](Sch-1) at (4.4,2.7) {\rotatebox{90}{$\copA$ comm}};
    \node[draw](Sch-88) at (-.6,-.5) {\rotatebox{90}{$\copA$ 2-TS}};
    \node[draw](Sch-408-654) at (2,5) {\rotatebox{90}{$\copA$ IP}};
    \node[draw](Sch-336) at (2,2.5) {\rotatebox{90}{$\copA$ EIP}};

    \node[draw](Sch-103) at (-3,8) {\rotatebox{90}{$\opB$ MCIP}};
    \node[draw](Sch-654) at (5,8) {\rotatebox{90}{$\opB$ MIP}};
    \node[draw](Sch-103-654) at (0,5) {\rotatebox{90}{$\opB$ 4-TS}};
    \node[draw](Sch-29) at (-3,5) {\rotatebox{90}{$\opB$ MICIP}};
    \node[draw](Sch-12) at (5,2.7) {\rotatebox{90}{$\opB$ comm}};
    \node[draw](Sch-18) at (0,-.5) {\rotatebox{90}{$\opB$ 2-TS}};
    \node[draw](Sch-408-432) at (2.6,5) {\rotatebox{90}{$\opB$ IP}};
    \node[draw](Sch-264) at (2.6,2.5) {\rotatebox{90}{$\opB$ EIP}};

    \node[draw](Sch-161) at (-2.4,8) {\rotatebox{90}{$\opC$ MCIP}};
    \node[draw](Sch-408) at (5.6,8) {\rotatebox{90}{$\opC$ MIP}};
    \node[draw](Sch-161-408) at (.6,5) {\rotatebox{90}{$\opC$ 4-TS}};
    \node[draw](Sch-183) at (-2.4,5) {\rotatebox{90}{$\opC$ MICIP}};
    \node[draw](Sch-6) at (5.6,2.7) {\rotatebox{90}{$\opC$ comm}};
    \node[draw](Sch-156) at (.6,-.5) {\rotatebox{90}{$\opC$ 2-TS}};
    \node[draw](Sch-432-654) at (3.2,5) {\rotatebox{90}{$\opC$ IP}};
    \node[draw](Sch-504) at (3.2,2.5) {\rotatebox{90}{$\opC$ EIP}};

    \node[draw](Sch-103-127) at (-5,4.5) {\rotatebox{90}{CIP}};
    \node[draw](Sch-23) at (-5,2.5) {\rotatebox{90}{ECIP}};
    \node[draw](Sch-29-93) at (-3,2.5) {\rotatebox{90}{ICIP}};
    \node[draw](Sch-5) at (-3,.5) {\rotatebox{90}{SS}};
    \node[draw](Sch-408-432-654) at (1.2,3) {\rotatebox{90}{all IP}};
    \node[draw](Sch-24) at (0,-2.5) {\rotatebox{90}{TS}};

    \draw[-{stealth}](Sch-432-654.north)--(Sch-654.south);
    \draw[-{stealth}](Sch-432-654.north)--(Sch-432.south);
    \draw[-{stealth}](Sch-432.north)--(Sch-);
    \draw[-{stealth}](Sch-654.north)--(Sch-);
    \draw[-{stealth}](Sch-408-432-654)--(Sch-432-654.south);
    \draw[-{stealth}](Sch-408-432-654)--(Sch-408-654.south);
    \draw[-{stealth}](Sch-408-432-654)--(Sch-408-432.south);
    \draw[-{stealth}](Sch-408-654.north)--(Sch-654.south);
    \draw[-{stealth}](Sch-408-654.north)--(Sch-408.south);
    \draw[-{stealth}](Sch-408-432.north)--(Sch-408.south);
    \draw[-{stealth}](Sch-408-432.north)--(Sch-432.south);
    \draw[-{stealth}](Sch-408.north)--(Sch-);
    \draw[-{stealth}](Sch-93.north)--(Sch-127.south);
    \draw[-{stealth}](Sch-127.north)--(Sch-);
    \draw[-{stealth}](Sch-103-127.north)--(Sch-103.south);
    \draw[-{stealth}](Sch-103-127.north)--(Sch-127.south);
    \draw[-{stealth}](Sch-103-127.north)--(Sch-161.south);
    \draw[-{stealth}](Sch-103.north)--(Sch-);
    \draw[-{stealth}](Sch-161.north)--(Sch-);
    \draw[-{stealth}](Sch-156.north)--(Sch-6.south);
    \draw[-{stealth}](Sch-156.north)--(Sch-408-432-654);
    \draw[-{stealth}](Sch-156.north)--(Sch-161-408.south);
    \draw[-{stealth}](Sch-156.north)--(Sch-183.south);
    \draw[-{stealth}](Sch-156.north)--(Sch-504.south);
    \draw[-{stealth}](Sch-504.north)--(Sch-432-654.south);
    \draw[-{stealth}](Sch-6.north)--(Sch-408.south);
    \draw[-{stealth}](Sch-161-408.north)--(Sch-408.south);
    \draw[-{stealth}](Sch-161-408.north)--(Sch-161.south);
    \draw[-{stealth}](Sch-183.north)--(Sch-161.south);
    \draw[-{stealth}](Sch-88.north)--(Sch-127-432.south);
    \draw[-{stealth}](Sch-88.north)--(Sch-408-432-654);
    \draw[-{stealth}](Sch-88.north)--(Sch-336.south);
    \draw[-{stealth}](Sch-88.north)--(Sch-1.south);
    \draw[-{stealth}](Sch-88.north)--(Sch-93.south);
    \draw[-{stealth}](Sch-336.north)--(Sch-408-654.south);
    \draw[-{stealth}](Sch-127-432.north)--(Sch-127.south);
    \draw[-{stealth}](Sch-127-432.north)--(Sch-432.south);
    \draw[-{stealth}](Sch-1.north)--(Sch-432.south);
    \draw[-{stealth}](Sch-264.north)--(Sch-408-432.south);
    \draw[-{stealth}](Sch-103-654.north)--(Sch-654.south);
    \draw[-{stealth}](Sch-103-654.north)--(Sch-103.south);
    \draw[-{stealth}](Sch-29.north)--(Sch-103.south);
    \draw[-{stealth}](Sch-18.north)--(Sch-12.south);
    \draw[-{stealth}](Sch-18.north)--(Sch-408-432-654);
    \draw[-{stealth}](Sch-18.north)--(Sch-103-654.south);
    \draw[-{stealth}](Sch-18.north)--(Sch-264.south);
    \draw[-{stealth}](Sch-18.north)--(Sch-29.south);
    \draw[-{stealth}](Sch-12.north)--(Sch-654.south);
    \draw[-{stealth}](Sch-23)--(Sch-103-127);
    \draw[-{stealth}](Sch-24.north)--(Sch-5.south);
    \draw[-{stealth}](Sch-24.north)--(Sch-156.south);
    \draw[-{stealth}](Sch-24.north)--(Sch-18.south);
    \draw[-{stealth}](Sch-24.north)--(Sch-88.south);
    \draw[-{stealth}](Sch-29-93)--(Sch-103-127);
    \draw[-{stealth}](Sch-29-93.north)--(Sch-183.south);
    \draw[-{stealth}](Sch-29-93.north)--(Sch-29.south);
    \draw[-{stealth}](Sch-29-93.north)--(Sch-93.south);
    \draw[-{stealth}](Sch-5)--(Sch-23);
    \draw[-{stealth}](Sch-5)--(Sch-29-93);
    \node at (-3.5,-1.95) {\footnotesize\begin{tabular}{@{}ll@{}}%
      SS&semi-symmetry\\%
      CIP&crossed inverse property\\%
      ECIP&equal CIP\\%
      ICIP&involutive CIP\\%
      MCIP&middle CIP\\%
      MICIP&middle involutive CIP\\%
    \end{tabular}};
    \node at (3.5,-1.6) {\footnotesize\begin{tabular}{@{}ll@{}}%
      MIP&middle inverse property\\%
      IP&inverse property\\%
      all IP&all inverse properties\\%
      EIP&equal IP\\\\%
      TS&total symmetry\\%
      2-TS&$\ZZ/2\ZZ$-twisted TS\\%
      4-TS&$\ZZ/4\ZZ$-twisted TS\\%
    \end{tabular}};
  \end{tikzpicture}}
  \caption{\label{fig:two-vars-lattice}Sublattice (of the full lattice of logically closed subsets of \num{990} laws) consisting of classes expressible in terms of two-variable equations.  The green region (middle) covers classes with conjectural pressure $p=5/6$.  The blue rectangle (left) $p=2/3$.  The orange rectangle (right) $p=1/2$.  Rather than left/middle/right (crossed) inverse properties I state middle properties for the appropriate operation, to obtain a more symmetric presentation.  The MCIP, MICIP, MIP, EIP, ECIP and commutation classes have a single outgoing edge hence cannot be written as conjunctions.}
  % python Schroeder_find_class.py hasse 1 5 6 12 18 23 24 29 88 93 103 127 156 161 183 264 336 408 432 504 654
\end{figure}

Note that quasigroup isomorphisms can be ignored when counting since it only introduces a factor of at most $n!<n^n$ so that it only adds a term $O(n\log n)$ in the exponential.
Since this conjecture is not central to this work, I postpone a discussion (and proof in some cases) to~\S\ref{sec:pressure}.

\subsection{Twists of commutativity (Class Lists \ref{class:empty} and~\ref{class:comm})}
\label{subsec:twist-comm}

We now arrive at the lists of classes.
See \autoref{fig:two-vars-lattice} for a depiction of the semi-lattice.
To state the upcoming properties, it is useful to introduce translation sets that appeared more abstractly in Proposition~\ref{prop:reduce}, as well as well-known inverse properties.

Given an element $x$ in a quasigroup~$Q$, there are six natural maps $Q\to Q$ that are pairwise inverses of each other:
\begin{itemize}
\item the left translations $y \mapsto x \opA y$ and $y \mapsto y \opC x$ (multiplication and division),
\item the middle translations $y \mapsto x \opC y$ and $y \mapsto y \opB x$,
\item the right translations $y \mapsto x \opB y$ and $y \mapsto y \opA x$ (division and multiplication).
\end{itemize}
The ``left'', ``middle'', and ``right'' translations are mapped cyclically by cyclic permutations of $\opA,\opB,\opC$.
Several notions then naturally come with a ``left'', ``middle'', and ``right'' versions.
The following result is an immediate corollary of Proposition~\ref{prop:reduce}, and is motivated by the point of view in~\cite{sokhatsky2020classification,sokhatsky2021classification}.

\begin{propdef}[Left, middle, right inverse properties]\label{propdef:lmr-inverse-prop}
The following four properties are equivalent, and the map $f:Q\to Q$ is automatically the same involution ($f\circ f=\mathrm{id}$) in each expression.
\begin{itemize}
\item[(i)] The law Sch-654 $x \opC (y \opB z) = (z \opB y) \opA x$ holds.
\item[(ii)] $(Q, \opA)$ has the \emph{left inverse property} (LIP): $f(x) \opA (x \opA y) = y$ for some~$f$.
\item[(iii)] $(Q, \opB)$ has the \emph{middle inverse property} (MIP): $f(x \opB y) = y \opB x$ for some~$f$.
\item[(iv)] $(Q, \opC)$ has the \emph{right inverse property} (RIP): $(y \opC x) \opC f(x) = y$ for some~$f$.
\item[(v)] One has $y \opC x = f(x) \opA y$ for some~$f$.
\item[(vi)] The sets of left multiplications $(x \opA \_)$ and left divisions $(\_ \opC x)$ coincide.
\end{itemize}
A quasigroup $(Q, \opA)$ obeying the left and right inverse properties $f(x) \opA (x \opA y) = y$ and $(y \opA x) \opA g(x) = y$ is said to have the \emph{inverse property} (IP).
If moreover the maps~$f,g$ appearing in the left/right inverse properties are equal then I shall dub it the \emph{equal inverse property} (EIP).
A quasigroup $(Q, \opA)$ obeying the left/middle/right inverse properties is said to have \emph{all inverse properties} (all IP).
\end{propdef}

\begin{remark}
In view of (iii), this property could also be dubbed $\ZZ/2\ZZ$-twisted commutativity of~$\opB$.
A quasigroup $(Q, \opB)$ obeying the middle inverse property with a trivial map $f=\mathrm{id}$ is simply a \emph{commutative} quasigroup.
\end{remark}

\begin{remark}
In loops, some of the distinctions between inverse properties collapse.  Evaluating the middle inverse property with a neutral element shows that $f=\mathrm{id}$ so the loop is commutative.  Evaluating the left \emph{or} right inverse properties with $y$ being the neutral element shows that the loop has two-sided inverses (given by the map~$f$).  The inverse property and equal inverse property are thus equivalent.  Exploring details of these equivalences would be interesting but it is out of scope of this work.  (See \cite{phillips2015scoop} for a discussion of inverse properties in unital magmas.)
\end{remark}

The weak inverse property (WIP), which states that $x \opA f(y \opA x) = f(y)$ for some map~$f$, is not captured by an equational law.  These properties have also been generalized further to $(r,s,t)$-inverse properties, in which different (composition) powers of~$f$ are inserted throughout the identity.

\begin{classlist}\label{class:empty}
A first class, with $p=0$.
\begin{itemize}
\item[\textbf{Sch-$\emptyset$}] (no eq)
There is not much to say about the variety $V(\emptyset)$ of \structure{all quasigroups}.  Denoted~$0$ in Schr\"oder's book~\cite{Schroder-1890}.
\end{itemize}
\end{classlist}

\begin{classlist}\label{class:comm}
Twists of commutativity (with pressure $p=1/2$ conjecturally).
\begin{itemize}

\item[\textbf{Sch-1}] ($\opA$ comm) $x \opA (y \opA z) = x \opA (z \opA y)$ iff $(Q, \opA)$ is \structure{commutative}
$\iff x \opA y \mkeqN{43} y \opA x
\iff x \opB y = y \opC x
\iff x \mkeqN{26} (x \opB y) \opB y
\iff x \mkeqN{16} y \opC (y \opC x)$.
Dubbed $C_1$ by Schr\"oder~\cite[p.636]{Schroder-1890}.

\item[Sch-6] ($\opC$ comm) $x \opA (y \opA z) = x \opA (z \opB y)$ iff $\opC$ is commutative (para Sch-1).
  Dubbed $C_2$ by Schr\"oder \cite[p.641]{Schroder-1890}.

\item[Sch-12] ($\opB$ comm) $x \opA (y \opA z) = x \opA (z \opC y)$ iff $\opB$ is commutative (para Sch-1).

\item[\textbf{Sch-408}] ($\opC$ MIP) $x \opA (y \opC z) = (z \opC y) \opB x$ iff $({\opB},{\opC},{\opA})$ have the \structure{(left, middle, right) inverse properties,} respectively,
in the sense that $f(x \opC y) = y \opC x$ for some map~$f$ (cf.\ Definition~\ref{propdef:lmr-inverse-prop}).  This map is automatically an involution, and $x \opB y = y \opA f(x)$ and $(x \opA y) \opA f(y) = x$.  Dubbed $J_2$ by Schr\"oder \cite[p.636]{Schroder-1890}.

\item[Sch-432] ($\opA$ MIP) $x \opB (y \opA z) = (z \opA y) \opC x$ iff $({\opC},{\opA},{\opB})$ have the (left, middle, right) inverse property, respectively (para Sch-408).  Dubbed $J_1$ by Schr\"oder~\cite[p.636]{Schroder-1890}.

\item[Sch-654] ($\opB$ MIP) $x \opC (y \opB z) = (z \opB y) \opA x$ iff (${\opA},{\opB},{\opC}$) have the (left, middle, right) inverse properties, respectively (para Sch-408).  Dubbed $J_3$ by Schr\"oder \cite[p.636]{Schroder-1890}.
\end{itemize}
\end{classlist}

\subsection{Twists of semi-symmetry (\autoref{class:ss})}
\label{subsec:twist-ss}

The following result is an immediate corollary of Proposition~\ref{prop:reduce}.

\begin{propdef}\label{propdef:lmr-crossed-inverse-prop}
The following four properties are equivalent, and the map~$f$ is automatically the same automorphism of~$Q$ in each expression.
\begin{itemize}
\item[(i)] The law Sch-103 $x \opA (y \opC z) = x \opB (y \opA z)$ holds.
\item[(ii)] $(Q, \opA)$ has the \emph{left crossed inverse property} (LCIP): $(x \opA y) \opA x = f(y)$ for some~$f$.
\item[(iii)] $(Q, \opB)$ has the \emph{middle crossed inverse property} (MCIP): $(x \opB y) \opB f(x) = y$ or equivalently $x \opB (y \opB f(x)) = y$ or equivalently $(f^{-1}(x) \opB y) \opB x = y$ or equivalently $f^{-1}(x) \opB (y \opB x) = y$ for some map~$f$.
\item[(iv)] $(Q, \opC)$ has the \emph{right crossed inverse property} (RCIP): $x \opC (y \opC x) = f^{-1}(y)$ for some~$f$.
\item[(v)] One has $x \opA y = x \opB f(y)$, or equivalently $x \opA y = f(x \opC y)$, or equivalently $f(x) \opC y = x \opB y$, for some~$f$.
\item[(vi)] The sets of middle translations $(x \opC \_)$ and right divisions $(x \opB \_)$ coincide.
\item[(vi)] The sets of middle translations $(\_ \opB x)$ and right multiplications $(\_ \opA x)$ coincide.
\end{itemize}
When the automorphism~$f$ is involutive ($f\circ f=\mathrm{id}$) I shall dub these properties \emph{(left, middle, right) involutive crossed inverse property} (LICIP, MICIP, RICIP).
\end{propdef}

\begin{remark}
To prove that $f$ is a morphism, use the LCIP three times in $f(x \opA y) = (f(y) \opA (x \opA y)) \opA f(y) = (((x \opA y) \opA x) \opA (x \opA y)) \opA f(y) = f(x) \opA f(y)$.
One also has $f(x \opC y) = x \opA y$.
In view of (v) the three operations are equal up to specific isotopies, so this property could also be dubbed $\ZZ$-twisted semi-symmetry.

Confusingly, the different characterizations of MCIP given here are sometimes called left and right~\cite{keedwell1999crossed}, which is somewhat justified in general magmas since they are no longer all equivalent~\cite{izbacs2007crossed}.
\end{remark}

\begin{propdef}\label{propdef:cip}
A quasigroup $Q$ obeying any two of the crossed inverse properties $(x\opA y)\opA x = f(y)$, $(x\opA y)\opA g(x) = y$, or $x\opA (h(y)\opA x) = y$ for some automorphisms $f,g,h$ obeys the third one and is said to have the \emph{crossed inverse property} (CIP).
The three automorphisms commute and $f(g(h(x)))=x$.

There exists a ring $(R,+,\times)$, equal to $\ZZ[X]/(1+X+X^2)$ or to one of its finite quotients, such that the map $\rho\colon R\to\Aut(Q)$,
\[
  a_0 + a_1 X + a_2 X^2 \mapsto f^{\circ a_0} \circ g^{\circ a_1} \circ h^{\circ a_2} ,
\]
defines a free action of the additive group $(R,+)$ on~$Q$ by quasigroup automorphisms.  The action is free in the sense that the only automorphism that has any fixed point is the image of $0\in R$.
\end{propdef}

\begin{proof}
By the criterion (vii) in Proposition~\ref{propdef:lmr-crossed-inverse-prop}, the three crossed inverse properties each state the equality of two of the three right-translation sets $\RT(\opA) = \{(\_\opA x)\}$ and $\RT(\opB) = \{(\_\opB x)\}$ and $\RT(\opC) = \{(\_\opC x)\}$, so any two of them implies the third equality.
The fact that $f,g,h$ commute is immediate since they are constructed from the quasigroup operations and are automorphisms: for instance, $g(f(x)) = g((x\opA x)\opA x) = (g(x)\opA g(x)) \opA g(x)=f(g(x))$.
The equality $f(g(h(x)))=x$ can be seen by
\begin{equation}
h^{-1}(y) \overset{\textnormal{RCIP}}= x \opA (y \opA x) \overset{\textnormal{MCIP}}= ((y\opA x) \opA g(y)) \opA (y \opA x) \overset{\textnormal{LCIP}}= f(g(y)) .
\end{equation}
Alternatively, using (v) in Proposition~\ref{propdef:lmr-crossed-inverse-prop}, one has
$f(g(h(x \opA y))) = f(g(x \opB y)) = f(x \opC y) = x \opA y$.

Because $f,g,h$ commute we can define a morphism $\tilde\rho$ from the additive group $\ZZ[X]$ to~$\Aut(Q)$ by mapping $X^{3k},X^{3k+1},X^{3k+2}$ to $f,g,h$, respectively.  The kernel $I=\ker\tilde\rho$ is an additive subgroup of $\ZZ[X]$, which we now show is a ring ideal.  It suffices to check that for any $P\in I$ one has $XP\in I$.
The key identities are
\begin{equation}\label{xfyeqgxy}
  \begin{aligned}
    x \opA f(y) & \overset{\textnormal{MCIP}}= g(x) \opB f(y) \overset{\textnormal{LCIP}}= g(x) \opA y, \\
    x \opA g(y) & \overset{\textnormal{RCIP}}= \bigl(y \opA (h(x) \opA y)\bigr) \opA g(y) \overset{\textnormal{MCIP}}= h(x) \opA y, \\
  \end{aligned}
\end{equation}
and their consequence
\begin{equation}
  x \opA h(y) = g(h(f(x))) \opA h(y) = h(f(x)) \opA f(h(y)) = f(x) \opA g(f(h(y))) = f(x) \opA y .
\end{equation}
Because of these three identities, we have
\begin{equation}\label{xtilderhoPy}
  x \opA \tilde\rho(P)(y) = \tilde\rho(XP)(x) \opA y .
\end{equation}
If $\tilde\rho(P)=\mathrm{id}$ (and more generally if it has any fixed point) then we can divide by~$y$ and learn that $\tilde\rho(XP)=\mathrm{id}$ too.
We are done proving that $I$ is a ring ideal, so $R=\ZZ[X]/I$ is a ring.

The ring~$R$ is a quotient of the ring of Einsenstein integers $R_{\textnormal{Eis}}=\ZZ[X]/(1+X+X^2)$ since $\tilde\rho(1+X+X^2)=f\circ g\circ h=\mathrm{id}$.  Such quotients (apart from the ring $R_{\textnormal{Eis}}$ itself) are all finite.  Indeed, consider an element of~$I$ that is not in the ideal $(1+X+X^2)$, and consider its residue $k+lX$ modulo $(1+X+X^2)$.  Since $R_{\textnormal{Eis}}$ is an integral domain, multiplication by $k+lX\neq 0$ is injective.  More precisely, this multiplication is a $\ZZ$-linear map on $R_{\textnormal{Eis}}$, with determinant $k^2-kl+l^2$, and the quotient $R_{\textnormal{Eis}}/(k+lX)$ has $k^2-kl+l^2$ elements.

The action of $R=\ZZ[X]/I$ is free since if $\tilde\rho(P)$ has a fixed point~$y_0$, then \eqref{xtilderhoPy} written for $y=y_0$ implies that $x=\tilde\rho(XP)(x)$ and then \eqref{xtilderhoPy} written for any other value of~$y$ implies $\tilde\rho(P)(y)=y$, namely $P\in I$.
\end{proof}

\begin{propdef}\label{propdef:icip-ecip}
In a crossed-inverse property (CIP) quasigroup~$Q$ (with the notation in Definition~\ref{propdef:cip}), if any of $f\circ f$, $g\circ g$, $h\circ h$ has a fixed point then $f$, $g$, $h$ are three involutions ($f\circ f=g\circ g=h\circ h=\mathrm{id}$) that describe an action of $(\ZZ/2\ZZ)^2$ on~$Q$.  I shall dub this the \emph{involutive crossed inverse property} (ICIP).
If any two of $f,g,h$ share a value, namely $f(x)=g(x)$ or $g(x)=h(x)$ or $h(x)=f(x)$ for some~$x$, then the functions $f=g=h$ are all equal to the same automorphism of order~$3$, and I shall dub this the \emph{equal crossed inverse property} (ECIP).
In an ECIP quasigroup~$Q$, either $f$ is the identity and the quasigroup is semi-symmetric, or $f$ has no fixed point and its orbits form a semi-symmetric quasigroup, so that $Q$ is an \emph{extension of a semi-symmetric quasigroup by $\ZZ/3\ZZ$} equipped with subtraction.
\end{propdef}

\begin{proof}
Recall from the proof of Proposition~\ref{propdef:cip} that mapping $1,X,X^2$ to $f,g,h$ defines a free action of $R=\ZZ[X]/I$ on~$Q$, where $I$ is a ring ideal containing the ideal $(1+X+X^2)$.  If $f\circ f$ or $g\circ g$ or $h\circ h$ has a fixed point then it must be the identity (the action is free), so that one of $2,2X,2X^2$ belongs to~$I$, hence $2 = 2X^3 + 2(1+X+X^2)(1-X)$ too.  Then $R$ is (a quotient of) $(\ZZ/2\ZZ)[X]/(1+X+X^2)$, which is a well-known description of the finite field with four elements~$\FF_4$.  For our purposes only the Abelian group structure $(\ZZ/2\ZZ)^2$ is important.  We thus have a free action of either $(\ZZ/2\ZZ)^2$ or $\ZZ/2\ZZ$ or the trivial group (in the case $f=g=h=\mathrm{id}$).

Regarding the ECIP property, note that $f(x)=g(x)$ is equivalent to $f^{-1}(g(x)) = x$, which implies that $X-1\in I$ (because the action is free), and likewise $g(x)=h(x)$ implies $X^2-X\in I$ and $f(x)=h(x)$ implies $1-X^2\in I$.  Multiplying by $1,X^2,X$, respectively, and reducing modulo $1+X+X^2$, these polynomials all reduce to $X-1$.  The ideal $(1+X+X^2,X-1)$ is equal to $(3,X-1)$, and the quotient of $\ZZ[X]$ by that ideal is then simply $\ZZ/3\ZZ$.  In detail, $f=g=h$ are automorphisms of order~$3$, as announced.
At this stage, $R$ is either $\ZZ/3\ZZ$ or the trivial ring.

If $f$ has a fixed point then (because the $R$-action is free) $R$ is the trivial ring, namely $f=g=h=\mathrm{id}$ and the quasigroup is semi-symmetric.  Otherwise, all orbits of~$f$ have three elements.  Let us show that $\opA$ respects the equivalence relation $x\sim x'$ if $x'=f^{\circ k}(x)$ for some $k\geq 0$.
The identity~\eqref{xfyeqgxy} reduces to $x \opA f(y) = f(x) \opA y$, from which one gets (using that $f$ is a morphism)
\begin{equation}
  x \opA f(y) = x \opA f(f(f(f(y)))) = f(f(x)) \opA f(f(y)) = f(f(x \opA y))
\end{equation}
and more generally for any $k,l\geq 0$,
\begin{equation}
  f^{\circ k}(x) \opA f^{\circ l}(y)
  = f^{\circ k}(x) \opA f^{\circ l+3(k+l)}(y)
  = f^{\circ 2(k+l)}(x) \opA f^{\circ 2(k+l)}(y)
  = f^{\circ 2(k+l)}(x \opA y) .
\end{equation}
Thus, for $x\sim x'$ and $y\sim y'$ one has $x\opA y\sim x'\opA y'$.  The quotient $(Q/\sim,\opA)$ is thus a quasigroup.  It is semi-symmetric since the twist by~$f$ disappears upon quotienting.
\end{proof}

\begin{classlist}\label{class:ss}
Twists of semi-symmetry (with $p=2/3$ conjecturally).
\begin{itemize}

\item[\textbf{Sch-5}] (SS) $x \opA (y \opA z) = x \opA (y \opB z)$ iff $Q$ is \structure{semi-symmetric}, that is, ${\opA}={\opB}={\opC}$, namely $x \mkeqN{14} y \opA (x \opA y)$ and the same for other operations.

\item[\textbf{Sch-23}] (ECIP) $x \opA (y \opA z) = x \opB (y \opB z)$ iff $Q$ has the \structure{equal crossed inverse property}, namely $B\colon x\mapsto x \opA (x \opA x)$ is a quasigroup automorphism of order~$3$ such that $x \opB y = B(x \opC y) = B(B(x \opA y))$ and $x \opA B(y) = B(x) \opA y$.  It has the crossed inverse property (see Proposition~\ref{propdef:cip}) with three equal automorphisms (see Proposition~\ref{propdef:icip-ecip}).  Equivalently, $x \mkeqN{692} y \opA (x \opA ((z \opA y) \opA z))$.  (Same equation for any of the six operations.)
Dubbed $C_{00}$ by Schr\"oder~\cite[p.639]{Schroder-1890}.

\item[\textbf{Sch-29}] ($\opB$ MICIP) $x \opA (y \opA z) = x \opB (y \opC z)$ iff $({\opA},{\opB},{\opC})$ have the \structure{(left, middle, right) involutive crossed inverse properties,} respectively (cf.\ Definition~\ref{propdef:lmr-crossed-inverse-prop}), in the sense that $x \opA y = x \opB f(y)$ and $x \opC y = f(x) \opB y$ for some involutive automorphism~$f$.  Equivalently, $x \mkeqN{2576} (y \opA ((z \opA x) \opA z)) \opA y$, which states that $f:x\mapsto (y\opA x)\opA y$ is $y$-independent and involutive.

\item[Sch-93] ($\opA$ MICIP) $x \opA (y \opB z) = x \opC (y \opC z)$ iff $({\opC},{\opA},{\opB})$ have the (left, middle, right) involutive crossed inverse properties, respectively (para Sch-29).

\item[Sch-183] ($\opC$ MICIP) $x \opB (y \opB z) = x \opC (y \opA z)$ iff $({\opB},{\opC},{\opA})$ have the (left, middle, right) involutive crossed inverse properties, respectively (para Sch-29).

\item[\textbf{Sch-103}] ($\opB$ MCIP) $x \opA (y \opC z) = x \opB (y \opA z)$ iff $({\opA},{\opB},{\opC})$ have the \structure{(left, middle, right) crossed inverse properties,} respectively (cf.\ Definition~\ref{propdef:lmr-crossed-inverse-prop}).  Equivalently $(x \opA x) \opA x \mkeqN{4588} (y \opA x) \opA y$.  Equivalently $x \opC (x \opC x) \mkeqN{4273} y \opC (x \opC y)$.  The map $x \mapsto (x \opA x) \opA x$ and its inverse $x \mapsto x \opC (x \opC x)$ are automorphisms.

\item[Sch-127] ($\opA$ MCIP) $x \opA (y \opC z) = x \opC (y \opB z)$ iff $({\opC},{\opA},{\opB})$ have the (left, middle, right) crossed inverse property (para Sch-103).
Dubbed $E_1$ by Schr\"oder~\cite[p.640]{Schroder-1890}.

\item[Sch-161] ($\opC$ MCIP) $x \opB (y \opA z) = x \opC (y \opB z)$ iff $({\opB},{\opC},{\opA})$ have the (left, middle, right) crossed inverse property (para Sch-103).

\item[\textbf{Sch-(29, 93)}] (ICIP) $x \opA (y \opA z) = x \opB (y \opC z)$ and $x \opA (y \opB z) = x \opC (y \opC z)$ iff $Q$ has the \structure{involutive crossed inverse property,} namely is semi-symmetric up to an action of $(\ZZ/2\ZZ)^2$ by automorphisms (see Proposition~\ref{propdef:icip-ecip}).

\item[\textbf{Sch-(103, 127)}] (CIP) $x \opA (y \opC z) = x \opB (y \opA z)$ and $x \opA (y \opC z) = x \opC (y \opB z)$ iff $Q$ has the \structure{crossed inverse property,} namely is semi-symmetric up to an action of Eisenstein integers (see Proposition~\ref{propdef:cip}).
\end{itemize}
\end{classlist}

\subsection{Twists of total symmetry (\autoref{class:ts})}
\label{subsec:twist-ts}

Finally, some conjunctions are strong enough to relate all six operations up to some unary maps.
\begin{classlist}\label{class:ts}
Twists of total symmetry (with $p=5/6$ conjecturally).
\begin{itemize}

\item[\textbf{Sch-24}] (TS) $x \opA (y \opA z) = x \opB (z \opB y)$ iff $Q$ is \structure{totally symmetric}, that is, ${\opA}={\opB}={\opC}$ are commutative, namely $x \mkeqN{492} y \opA (x \opA (z \opA (z \opA y)))$ and the same for other operations.
Dubbed $Q_0$ by Schr\"oder~\cite[p.644]{Schroder-1890}.

\item[\textbf{Sch-18}] ($\opB$ 2-TS) $x \opA (y \opA z) = x \opB (z \opA y)$ iff $(Q,\opB)$ is \structure{$\ZZ/2\ZZ$-twisted totally symmetric} in the sense that all six operations are related by $x \opA y = y \opC x = x \opB f(y)$ and $x\opB y=y\opB x$ and $f(x \opA y) = y \opA x$ and $f(x \opC y) = y \opC x$ for some involutive automorphism~$f$.  Equivalently, $(Q,\opB)$ is commutative (Sch-12) and has the inverse property (Sch-408, Sch-432).

\item[Sch-88] ($\opA$ 2-TS) $x \opA (y \opB z) = x \opC (z \opB y)$ iff $(Q,\opA)$ is $\ZZ/2\ZZ$-twisted totally symmetric (para Sch-18).
Dubbed $Q_1$ by Schr\"oder~\cite[p.644]{Schroder-1890}.

\item[Sch-156] ($\opC$ 2-TS) $x \opB (y \opA z) = x \opC (z \opA y)$ iff $(Q,\opC)$ is $\ZZ/2\ZZ$-twisted totally symmetric (para Sch-18).

\item[\textbf{Sch-264}] ($\opB$ EIP) $x \opA (y \opA z) = (z \opA y) \opB x$ iff $(Q, \opB)$ has the \structure{equal inverse property,} namely has the left and right inverse properties (cf.\ Definition~\ref{propdef:lmr-inverse-prop}) with the same map~$f$.  Equivalently, $f(x \opA y) = y \opA x$ and $f(x \opC y) = y \opC x$.  The map~$f$ is automatically an involution and an anti-automorphism of $(Q, \opB)$, that is, $f(x \opB y) = f(y) \opB f(x)$.  All six operations are related by the above identities and $x \opA y = f(y) \opB x$ and $x \opC y = y \opB f(x)$, but this is weaker than being $\ZZ/2\ZZ$-twisted totally symmetric (Sch-18).

\item[Sch-336] ($\opA$ EIP) $x \opA (y \opB z) = (z \opB y) \opB x$ iff $(Q, \opA)$ has the equal inverse property (para Sch-264).  Dubbed $K_1$ by Schr\"oder~\cite[p.636]{Schroder-1890}.

\item[Sch-504] ($\opC$ EIP) $x \opB (y \opB z) = (z \opB y) \opC x$ iff $(Q, \opC)$ has the equal inverse property (para Sch-264).

\item[\textbf{Sch-(103, 654)}] ($\opB$ 4-TS) $x \opA (y \opC z) = x \opB (y \opA z)$ and $x \opC (y \opB z) = (z \opB y) \opA x$ iff $(Q,\opB)$ is \structure{$\ZZ/4\ZZ$-twisted totally symmetric} in the sense that $C\colon x\mapsto (x \opA x) \opA x$ is an automorphism of order~$4$ and all six operations are related by $x \opB y = C(x) \opC y = x \opA C^{-1}(y)$ and $x \opB y = C(C(y \opB x))$.  This reduces to Sch-18 when $C\circ C=\mathrm{id}$.

\item[Sch-(127, 432)] ($\opA$ 4-TS) $x \opA (y \opC z) = x \opC (y \opB z)$ and $x \opB (y \opA z) = (z \opA y) \opC x$ iff $(Q,\opA)$ is $\ZZ/4\ZZ$-twisted totally symmetric (para Sch-(103, 654)).

\item[Sch-(161, 408)] ($\opC$ 4-TS) $x \opB (y \opA z) = x \opC (y \opB z)$ and $x \opA (y \opC z) = (z \opC y) \opB x$ iff $(Q,\opC)$ is $\ZZ/4\ZZ$-twisted totally symmetric (para Sch-(103, 654)).

\item[\textbf{Sch-(408, 432)}] ($\opB$ IP) $x \opA (y \opC z) = (z \opC y) \opB x$ and $x \opB (y \opA z) = (z \opA y) \opC x$ iff $(Q, \opB)$ has the \structure{inverse property,} namely has the left and right inverse properties (cf.\ Definition~\ref{propdef:lmr-inverse-prop}).  Equivalently, $f(x \opA y) = y \opA x$ and $h(x \opC y) = y \opC x$ for some maps~$f,h$ which are automatically involutions.
The six operations are related by these identities and by $x \opB y = f(h(x) \opA y)$ and $x \opC y = f(x) \opA h(y)$ and $y \opB x = f(h(x) \opB h(y))$.  If the two involutions coincide one has the equal inverse property (Sch-264).

\item[Sch-(408, 654)] ($\opA$ IP) $x \opA (y \opC z) = (z \opC y) \opB x$ and $x \opC (y \opB z) = (z \opB y) \opA x$ iff $(Q, \opA)$ has the inverse property (para Sch-(408, 432)).

\item[Sch-(432, 654)] ($\opC$ IP) $x \opB (y \opA z) = (z \opA y) \opC x$ and $x \opC (y \opB z) = (z \opB y) \opA x$ iff $(Q, \opC)$ has the inverse property (para Sch-(408, 432)).

\item[\textbf{Sch-(408, 432, 654)}] (all IP) $x \opA (y \opC z) = (z \opC y) \opB x$ and $x \opB (y \opA z) = (z \opA y) \opC x$ and $x \opC (y \opB z) = (z \opB y) \opA x$ states that all three operations have the left, middle, and right inverse properties.  These involve three involutive automorphisms, defined by $f(x \opA y) = y \opA x$ and $g(x \opB y) = y \opB x$ and $h(x \opC y) = y \opC x$, which commute and obey $h = f \circ g$.  Collectively, these describe an action of $(\ZZ/2\ZZ)^2$ on the quasigroup~$Q$ by automorphisms.  An example where all three are pairwise distinct is given in~\eqref{quasigroup-all-IP}.  One has $x \opB y = g(x) \opA f(y)$ and $x \opC y = f(x) \opA h(y)$.
\end{itemize}
\end{classlist}

\section{Groups (3 classes) (\autoref{class:gp})}
\label{sec:groups}

Next comes a well-known class (and its parastrophes), which turns out to be an outlier.  Indeed, all other classes are either weaker than total symmetry (treated in \S\ref{sec:twovars}), or imply an Abelian group structure on~$Q$ up to a choice of origin, in other words an affine $\ZZ$-module structure on~$Q$.
The only classes that imply associativity are Sch-3 (Abelian groups) and Sch-7 (Boolean groups).
The only classes implied by associativity of~$\opA$ are the equal inverse property Sch-336 ($\opA$~EIP) and its consequences, which can be read in \autoref{fig:two-vars-lattice}, namely the inverse property Sch-(408, 654) ($\opA$~IP), and each of the middle inverse properties Sch-408 ($\opC$~MIP) and Sch-654 ($\opB$~MIP), as well as the no-equation class Sch-.
The sublattice generated by associativity and its parastrophes is given in \autoref{fig:sublattice-associativity} on \autopageref{fig:sublattice-associativity}.
\begin{classlist}\label{class:gp}
Associativity and its parastrophes.
\begin{itemize}

\item[\textbf{Sch-65}] ($\opA$ Grp) $x \opA (y \opB z) = y \opB (x \opA z)$ iff $(Q,\opA)$ is a \structure{group}, namely $\opA$ is associative;
  then $x \opB y = y x^{-1}$ and $x \opC y = y^{-1}x$.
  Dubbed $A_1$ by Schr\"oder \cite[p.635]{Schroder-1890}.

\item[Sch-105] ($\opB$ Grp) $x \opA (y \opC z) = y \opB (x \opA z)$ iff $(Q, \opB)$ is a group (para Sch-65).

\item[Sch-123] ($\opC$ Grp) $x \opA (y \opC z) = y \opC (x \opA z)$ iff $(Q, \opC)$ is a group (para Sch-65).
  Dubbed $A_2$ by Schr\"oder \cite[p.641]{Schroder-1890}.
\end{itemize}
\end{classlist}

\section{Isotopes of Abelian groups}
\label{sec:isotopy}

\subsection{Proof of the Four Quasigroups Theorem}
\label{subsec:prove-four-quasigroups}

All remaining classes (to be listed in \S\ref{subsec:quasilinear}--\S\ref{subsec:medial}) are described by generalized associativity laws of the form
\begin{equation}\label{proof-gen-assoc}
  x \odot_1 (y \odot_2 z) = (x \odot_3 y) \odot_4 z
\end{equation}
for $\odot_i$~among the six quasigroup operations ${\opA},{\opB},{\opC},{\tilde\opA},{/},{\backslash}$.
We begin here with a proof of the four quasigroups theorem, already discussed in~\S\ref{subsec:four-quasigroups}, which introduces a group operation $(Q,\times)$.  Then in \S\ref{subsec:Belousov-criterion} we explain Belousov's criterion and in \S\ref{subsec:underlying-commutativity} how it forces $(Q,\times)$ to be Abelian.
We then list in \S\ref{subsec:quasilinear} the classes that do not force the quasigroup operations to be linear, postponing the non-commutative and commutative linear models to \S\ref{subsec:noncomm} and~\S\ref{subsec:medial}, respectively.

\begin{prop}[Four quasigroups theorem]\label{prop:four-quasigroups}
  If four quasigroups $(Q,\odot_i)$ on the same carrier~$Q$ satisfy~\eqref{proof-gen-assoc} then the four operations are isotopes of the same group $(Q,\times)$, and specifically $x \odot_1 y = \alpha(x) \times \mu(y)$, and $x \odot_2 y = \mu^{-1}(\beta(x) \times \gamma(y))$, and $x \odot_3 y = \nu^{-1}(\alpha(x) \times \beta(y))$, and $x \odot_4 y = \nu(x) \times \gamma(y)$ for some bijections $\alpha,\beta,\gamma,\mu,\nu$.
\end{prop}

\begin{proof}
Denote by $\oslash_i$ and $\obslash_i$ the left and right divisions associated to $\odot_i$, namely $(x \odot_i y) \oslash_i y = x$ and $x \obslash_i (x \odot_i y) = y$ (note that we do not use conjugate divisions here for notational convenience).

We fix an element $w\in Q$ once and for all and define several bijections of~$Q$,
\begin{equation}
  \begin{aligned}
    \mu(x) & = w \odot_1 x , \qquad
    & \alpha(x) & = x \odot_1 (w \odot_2 w) = (x \odot_3 w) \odot_4 w ,
    \\
    && \beta(x) & = w \odot_1 (x \odot_2 w) = (w \odot_1 x) \odot_2 w ,
    \\
    \nu(x) & = x \odot_4 w , \qquad
    & \gamma(x) & = w \odot_1 (w \odot_2 x) = (w \odot_3 w) \odot_4 x ,
  \end{aligned}
\end{equation}
and a binary operation
\begin{equation}
  x \times y = \alpha^{-1}(x) \odot_1 \mu^{-1}(y) .
\end{equation}
The specialization of~\eqref{proof-gen-assoc} to $x=w$ relates $\odot_2$ to~$\odot_4$ up to one-variable functions, the specialization to $y=w$ relates $\odot_4$ to~$\odot_1$, and the specialization to $z=w$ relates $\odot_3$ to~$\odot_1$.  In terms of the operation~$\times$, this yields
\begin{equation}
  \begin{aligned}
    x \odot_1 y & = \alpha(x) \times \mu(y) , \\
    x \odot_2 y & = w \obslash_1 (\alpha((w \odot_3 x) \oslash_3 w) \times \mu(w \odot_2 y)) = \mu^{-1}(\beta(x) \times \gamma(y)) , \\
    x \odot_3 y & = (\alpha(x) \times \mu(y \odot_2 w)) \oslash_4 w = \nu^{-1}(\alpha(x) \times \beta(y)) , \\
    x \odot_4 y & = \alpha(x \oslash_3 w) \times \mu(w \odot_2 y) = \nu(x) \times \gamma(y),
  \end{aligned}
\end{equation}
which coincides with~\eqref{four-quasigroups-result}.
There remains to prove that $\times$ is associative.

We evaluate~\eqref{proof-gen-assoc} with the expressions of all four operations to get
\begin{equation}
  \alpha(x) \times \mu(\mu^{-1}(\beta(y) \times\gamma(z)))
  = \nu(\nu^{-1}(\alpha(x)\times\beta(y))) \times \gamma(z) .
\end{equation}
After changing variables to $x'=\alpha(x)$, $y'=\beta(y)$, $z'=\gamma(z)$, this is precisely the usual associativity condition $x'\times(y'\times z')=(x'\times y')\times z'$.  To finish, we observe that $\times$ is an isotope of a quasigroup operation hence it is a quasigroup operation.  Thus $(Q,\times)$ is a group as announced.
\end{proof}

Observe that $\alpha(w) = \beta(w) = \gamma(w) = w \odot_1 (w \odot_2 w) = (w \odot_3 w) \odot_4 w$ and that this is the neutral element of $(Q,\times)$ since
\begin{equation}
  \alpha(w) \times x = \alpha^{-1}(\alpha(w)) \odot_1 \mu^{-1}(x)
  = w \odot_1 \mu^{-1}(x) = \mu(\mu^{-1}(x)) = x .
\end{equation}
The group operation and neutral element depend on the arbitrary choice of $w\in Q$ that we made initially.

\subsection{The Belousov criterion}
\label{subsec:Belousov-criterion}

The four quasigroups theorem on generalized associativity is conveniently completed by a result on commutativity, and more generally a description of when two group isotopes are equal.  Proposition~\ref{prop:Belousov} is a special case of Proposition~\ref{prop:isotope-groups} below (whose special cases rely on Lemma~\ref{lem:inverse-automorphism}).  After proving these results, we explain how they apply in our setting.

\begin{prop}[Belousov's criterion]\label{prop:Belousov}
  If a group $(Q,\times)$ is such that there exist four bijections $\alpha,\beta,\gamma,\delta$ with $\alpha(x)\times\beta(y) = \gamma(y)\times\delta(x)$ for all $x,y$, then the group is Abelian.
\end{prop}

\begin{prop}[Group isotopy]\label{prop:isotope-groups}
  Consider a group $(Q,\times)$ with neutral element~$e$, and six bijections $\alpha,\beta,\gamma,\delta,\mu,\nu$.
  Assume that for all $x,y$,
  \[
    \mu(\alpha(x) \times \beta(y)) = \nu(\gamma(x) \times \delta(y)) .
  \]
  Then there exists a group automorphism $\theta$ of $(Q,\times)$ and elements $a,b\in Q$ such that
  \[
    \alpha(x) = a \times \theta(\gamma(x)) \times b^{-1} , \quad
    \beta(x) = b \times \theta(\delta(x)) , \quad
    \mu(a \times \theta(x)) = \nu(x) .
  \]
  Furthermore, the following six special cases occur, where $\iota\colon x\mapsto x^{-1}$ is the inversion map in $(Q,\times)$:
  \begin{itemize}
  \item if $\alpha=\gamma$, then $a=b$ and $\theta(x) = a^{-1}\times x\times a$;
  \item if $\beta=\delta$, then $b=e$ and $\theta=\mathrm{id}$;
  \item if $\mu=\nu$, then $a=e$ and $\theta=\mathrm{id}$;
  \item if $\alpha=\iota\circ\gamma$ then $b=a$ and $\theta(x) = x^{-1}$ and the group is Abelian;
  \item if $\beta=\iota\circ\delta$ then $b=e$ and $\theta(x) = x^{-1}$ and the group is Abelian;
  \item if $\mu=\nu\circ\iota$ then $a=e$ and $\theta(x) = x^{-1}$ and the group is Abelian.
  \end{itemize}
\end{prop}

\begin{lemma}\label{lem:inverse-automorphism}
  If for some $a,b\in Q$ the map $\theta\colon x\mapsto a \times x^{-1} \times b^{-1}$ is a group (auto)morphism then $(Q,\times)$ is Abelian, $a=b$ and $\theta(x)=x^{-1}$.
\end{lemma}

\begin{proof}
  Group morphisms preserve the identity element~$e$ so $a\times b^{-1} = e$, namely $b=a$.  Then the morphism property implies commutativity since $x \times y = \theta^{-1}(\theta(x) \times \theta(y)) = \theta^{-1} (a \times x^{-1} \times y^{-1} \times a^{-1}) = y \times x$.
\end{proof}

\begin{proof}[Proof of Proposition~\ref{prop:isotope-groups}]
  By changing variables as $x\to\gamma(x)$ and $y\to\delta(y)$ the identity is put in the form $A(x) \times B(y) = C(x \times y)$ with $A=\alpha\circ\gamma^{-1}$, $B=\beta\circ\delta^{-1}$, and $C=\mu^{-1}\circ\nu$.  Specializing the identity to $y=e$ and $x=e$, where $e$ denotes the group identity in $(Q,\times)$, one gets $A(x) = C(x) B(e)^{-1}$ and $B(x) = A(e)^{-1} C(x)$.  Introduce the elements $a,b\in Q$ and bijection~$\theta$
  \begin{equation}
    a = C(e) , \qquad b = B(e) , \qquad \theta(x) = C(e)^{-1} \times C(x) .
  \end{equation}
  In terms of these objects (and using that $A(e) B(e) = C(e)$), one has
  \begin{equation}
    A(x) = a \times \theta(x) \times b^{-1} , \qquad
    B(x) = b \times \theta(x) , \qquad
    C(x) = a \times \theta(x) ,
  \end{equation}
  and the identity we started with reads
  \begin{equation}
    (a \times \theta(x) \times b^{-1}) \times (b \times \theta(y)) = a \times \theta(x \times y) ,
  \end{equation}
  which after cancellations yields $\theta(x) \times \theta(y) = \theta(x \times y)$, namely $\theta$ is a group homomorphism.

  If $\alpha=\gamma$ then $x = a\times\theta(x)\times b^{-1}$.  Specializing to $x=e$ this reduces to $a\times b^{-1} = e$ so $a=b$, and the original identity immediately gives $\theta(x) = a^{-1}\times x \times a$.  If $\beta=\delta$ then $x = b\times\theta(x)$, which by the same reasoning gives $b=e$ and $\theta(x) = x$.  The logic is identical for the case $\mu=\nu$.

  We turn to the fourth special case.  If $\alpha=\iota\circ\gamma$ then $x = a \times \theta(x^{-1}) \times b^{-1}$, thus $\theta(x) = a^{-1}\times x^{-1}\times b$, to which Lemma~\ref{lem:inverse-automorphism} applies.
  The other two special cases are identical with $(a,b)$ replaced by $(b,e)$ and $(a,e)$, respectively.
\end{proof}

\begin{proof}[Proof of Proposition~\ref{prop:Belousov}]
  The identity $\alpha(x)\times\beta(y) = \gamma(y)\times\delta(x)$ can be rewritten $\mathrm{id}(\alpha(x)\times\beta(y)) = \iota((\iota\circ\delta)(x) \times (\iota\circ\gamma)(y))$, which fits in the last special case of Proposition~\ref{prop:isotope-groups}.  Thus, $(Q,\times)$ is Abelian and there exists $b\in Q$ such that
  \begin{equation}
    \alpha(x) = \delta(x) \times b^{-1} , \quad
    \beta(x) = b \times \gamma(x) . \qedhere
  \end{equation}
\end{proof}

Consider now the generalized associativity~\eqref{proof-gen-assoc} in the case where the four operations~$\odot_i$ are among the six quasigroup operations of a single quasigroup $(Q,\opA)$.  The four quasigroups theorem describes the four quasigroups $(Q,\odot_i)$ as isotopes of a group $(Q,\times)$, and thus describe $(Q,\opA)$ in four different ways as isotopes of $(Q,\times)$.  Specifically,
\begin{equation}
  x \opA y = \xi_i^{-1} (\zeta_i(x) \times \eta_i(y)) , \qquad i=1,2,3,4 ,
\end{equation}
where the bijections $\xi_i$, $\zeta_i$, $\eta_i$ are related to $\alpha,\beta,\gamma,\delta,\mu,\nu$ according to Table~\ref{tab:xiizetaietai} that lists, for each~$i$, and each of the six possibilities for~$\odot_i$, what $\xi_i,\zeta_i,\eta_i$ are.
Here, $\iota$~denotes the inversion map in $(Q,\times)$.

\begin{table}\centering
\caption{\label{tab:xiizetaietai}Bijections such that $\xi_i(x \opA y) = \zeta_i(x) \times \eta_i(y)$ according to which of the six quasigroup equations the operation~$\odot_i$ coincides with, for $i=1,2,3,4$.}
\begin{tabular}{>$c<$>$c<$}
\begin{array}{lccc}
  \toprule
  & \xi_1 & \zeta_1 & \eta_1 \\ \midrule
  \odot_1 = \opA & \mathrm{id} & \alpha & \mu \\
  \odot_1 = \opB & \iota\circ\mu & \iota & \alpha \\
  \odot_1 = \opC & \iota\circ\alpha & \mu & \iota \\
  \odot_1 = \tilde\opA & \iota & \iota\circ\mu & \iota\circ\alpha \\
  \odot_1 = / & \alpha & \mathrm{id} & \iota\circ\mu \\
  \odot_1 = \backslash & \mu & \iota\circ\alpha & \mathrm{id}
  \\\bottomrule
\end{array}
&
\begin{array}{lccc}
  \toprule
  & \xi_2 & \zeta_2 & \eta_2 \\ \midrule
  \odot_2 = \opA & \mu & \beta & \gamma \\
  \odot_2 = \opB & \iota\circ\gamma & \iota\circ\mu & \beta \\
  \odot_2 = \opC & \iota\circ\beta & \gamma & \iota\circ\mu \\
  \odot_2 = \tilde\opA & \iota\circ\mu & \iota\circ\gamma & \iota\circ\beta \\
  \odot_2 = / & \beta & \mu & \iota\circ\gamma \\
  \odot_2 = \backslash & \gamma & \iota\circ\beta & \mu
  \\\bottomrule
\end{array}
\\[-.8pt]
\begin{array}{lccc}
  \toprule
  & \xi_3 & \zeta_3 & \eta_3 \\ \midrule
  \odot_3 = \opA & \nu & \alpha & \beta \\
  \odot_3 = \opB & \iota\circ\beta & \iota\circ\nu & \alpha \\
  \odot_3 = \opC & \iota\circ\alpha & \beta & \iota\circ\nu \\
  \odot_3 = \tilde\opA & \iota\circ\nu & \iota\circ\beta & \iota\circ\alpha \\
  \odot_3 = / & \alpha & \nu & \iota\circ\beta \\
  \odot_3 = \backslash & \beta & \iota\circ\alpha & \nu
  \\\bottomrule
\end{array}
&
\begin{array}{lccc}
  \toprule
  & \xi_4 & \zeta_4 & \eta_4 \\ \midrule
  \odot_4 = \opA & \mathrm{id} & \nu & \gamma \\
  \odot_4 = \opB & \iota\circ\gamma & \iota & \nu \\
  \odot_4 = \opC & \iota\circ\nu & \gamma & \iota \\
  \odot_4 = \tilde\opA & \iota & \iota\circ\gamma & \iota\circ\nu \\
  \odot_4 = / & \nu & \mathrm{id} & \iota\circ\gamma \\
  \odot_4 = \backslash & \gamma & \iota\circ\nu & \mathrm{id}
  \\\bottomrule
\end{array}
\end{tabular}
\end{table}

For instance, this table can be used for Sch-54 $x \opA (y \opB z) = y \opA (x \opB z)$ by writing it (after renaming the variables $y,z$) as $x \opA (y / z) = (x \opB y) \tilde\opA z$, from which we see $\odot_1=\opA$ and $\odot_2=/$ and $\odot_3=\opB$ and $\odot_4=\tilde\opA$ hence
\begin{equation}
\begin{aligned}
x \opA y & = \alpha(x) \times \mu(y)
= \beta^{-1}(\mu(x) \times \gamma(y)^{-1})
\\
& = \beta^{-1}(\alpha(y)^{-1} \times \nu(x))
= \nu(y) \times \gamma(x) .
\end{aligned}
\end{equation}
From $\alpha(x) \times \mu(y) = \nu(y) \times \gamma(x)$, or from $\beta^{-1}(\mu(x) \times \gamma(y)^{-1}) = \beta^{-1}(\alpha(y)^{-1} \times \nu(x))$, Belousov's criterion implies that $(Q,\times)$ is Abelian.

\subsection{Commutativity of the underlying group}
\label{subsec:underlying-commutativity}

The four quasigroups theorem and Belousov's criterion apply to all remaining classes, which means that $(Q,\times)$ is always Abelian.  In the following proposition we keep multiplicative notation, but in subsequent sections we switch to additive notation $(Q,+)$ instead to emphasize commutativity.

\begin{prop}[Trichotomy of equivalence classes]\label{prop:trichotomy}
  All equivalence classes of conjunctions of the \num{990} Schr\"oder laws are either consequences of total symmetry (and listed in \S\ref{sec:twovars}), or equivalent to a parastrophe of associativity (listed in \S\ref{sec:groups}), or imply that the quasigroup operation takes the form $x\opA y = \xi^{-1}(\varphi(x) \times \psi(y))$ for some bijections $\varphi,\psi,\xi:Q\to Q$, where $(Q,\times)$ is an Abelian group.
\end{prop}

\begin{proof}
By inspection (a general argument is desirable), each of the classes below includes at least one equation in the following list, so that $(Q,\times)$ is Abelian in all cases.
\begin{itemize}
\item Any equation of the form $x \opA (\_ ? \_) = \_ \opA (\_ ? \_)$ has $\odot_1=\opA$ and $\odot_4=\tilde\opA$ hence $x \opA y = \alpha(x) \times \mu(y) = \nu(y) \times \gamma(x)$, for which Belousov's criterion makes $(Q,\times)$ Abelian.

\item By parastrophy, any equation of the form $x \opB (\_ ? \_) = \_ \opB (\_ ? \_)$ or $x \opC (\_ ? \_) = \_ \opC (\_ ? \_)$ also implies that $(Q,\times)$ is Abelian.  This can also be seen from the tables: $(\odot_1,\odot_4) = (\opB,/)$ or $(\opC,\backslash)$ implies $(\zeta_1,\zeta_4)=(\iota,\mathrm{id})$ or $(\eta_1,\eta_4)=(\iota,\mathrm{id})$, which fit in the special cases of Proposition~\ref{prop:isotope-groups}.

\item For Sch-20, $x \opA (y \opA z) = y \opB (z \opA x)$, one has $\odot_3=\tilde\opA$ and $\odot_4=/$ hence $\xi_3=\iota\circ\nu$ and $\xi_4=\nu$, which implies that $(Q,\times)$ is Abelian by the last special case in Proposition~\ref{prop:isotope-groups}.

\item For Sch-34 $x \opA (y \opA z) = z \opB (y \opC x)$, one has $\odot_1 = \opA$ and $\odot_3 = \backslash$ hence $\zeta_1=\alpha$ and $\zeta_3=\iota\circ\alpha$, which fits in the fourth special case in Proposition~\ref{prop:isotope-groups}.

\item For Sch-25, $x \opA (y \opA z) = y \opB (x \opB z)$, one has $\odot_1=\opA$ and $\odot_2=\tilde\opA$ and $\odot_3=\opB$ and $\odot_4=/$ hence
\[
\begin{aligned}
x \opA y & = \alpha(x) \times \mu(y)
= \mu^{-1}(\beta(y) \times \gamma(x))
\\
& = \beta^{-1}(\alpha(y)^{-1} \times \nu(x))
= \nu^{-1}(x \times \gamma(y)^{-1}) .
\end{aligned}
\]
Proposition~\ref{prop:isotope-groups} applied to $\alpha(x) \times \mu(y) = \nu^{-1}(x \times \gamma(y)^{-1})$ and to $\beta^{-1}(\alpha(y)^{-1} \times \nu(x)) = \nu^{-1}(x \times \gamma(y)^{-1})$ yield the existence of elements $a,b\in Q$ and automorphisms $\theta_1,\theta_2$ such that (among other properties)
\[
a \times \theta_1(x) = \nu^{-1}(x) , \qquad
\nu(x) = b \times \theta_2(x^{-1}) ,
\]
from which one gets $x = a \times \theta_1(\nu(x)) = a \times \theta_1(b) \times \theta_1(\theta_2(x^{-1}))$, so that the automorphism $\theta_1\circ\theta_2$ is given by $\theta_1(\theta_2(x)) = \theta_1(b)^{-1} \times a^{-1} \times x^{-1}$.  Lemma~\ref{lem:inverse-automorphism} then concludes that $(Q,\times)$ is Abelian.

\item For Sch-38 $x \opA (y \opA z) = y \opC (z \opA x)$ the logic is similar: the consequences of $\odot_1=\opA$ and $\odot_2=\tilde\opA$ give two expressions of $x\opA y$, and their equality falls in the purview of Proposition~\ref{prop:isotope-groups}, which relates in particular $\mu$ to a first automorphism up to the inversion~$\iota$, and the consequences of $\odot_1=\opA$ and $\odot_4=\backslash$ give another identification of $\mu$ with an automorphism without the inversion~$\iota$:
\[
\begin{aligned}
x \opA y = \alpha(x) \times \mu(y) = \mu^{-1}(\beta(y) \times \gamma(x)) & \implies a \times \theta_1(x^{-1}) = \mu^{-1}(x) ,
\\
x \opA y = \alpha(x) \times \mu(y) = \gamma(\nu(x)^{-1} \times y) & \implies \mu(x) = b \times \theta_2(x) ,
\end{aligned}
\]
hence $x = a \times \theta_1((b\times\theta_2(x))^{-1}) = a \times \theta_1(\theta_2(x^{-1})) \times \theta_1(b)^{-1}$.  Once more Lemma~\ref{lem:inverse-automorphism} applies to the automorphism $\theta_1\circ\theta_2$ and we learn that $(Q,\times)$ is Abelian.

\item For Sch-79 $x \opA (y \opB z) = z \opB (x \opC y)$ one proceeds identically: one has
\[
\begin{array}{lccc}
  & \xi_i & \zeta_i & \eta_i \\ \midrule
  \odot_1 = \opA & \mathrm{id} & \alpha & \mu \\
  \odot_2 = \opB & \iota\circ\gamma & \iota\circ\mu & \beta \\
  \odot_3 = \opC & \iota\circ\alpha & \beta & \iota\circ\nu \\
  \odot_4 = / & \nu & \mathrm{id} & \iota\circ\gamma \\
\end{array}
\]
Proposition~\ref{prop:isotope-groups} identifies all entries in a given column to each other up to automorphisms and left/right multiplications.  In particular, the $\xi_i$ column features two entries $\mathrm{id}$ and $\iota\circ\alpha$ and the $\zeta_i$ column features $\alpha$ and $\mathrm{id}$, so in the end $\alpha$ and $\iota\circ\alpha$ are identified up to an automorphism and left/right multiplication.  Lemma~\ref{lem:inverse-automorphism} then concludes that $(Q,\times)$ is Abelian.

\end{itemize}
\end{proof}

\subsection{Non-linear Abelian group isotopes (13 classes) (Class Lists~\ref{class:nonlinear-abelian} and~\ref{class:nonlinear-boolean})}
\label{subsec:quasilinear}

As summarized in Proposition~\ref{prop:trichotomy}, there only remains to describe varieties of quasigroups on which one can define an operation~$+$ such that $(Q,+)$ is an Abelian group and $\opA$ is an isotope of~$+$, namely $x \opA y = \gamma^{-1}(\alpha(x)+\beta(y))$ for some bijections $\alpha,\beta,\gamma:Q\to Q$ (renamed compared to Proposition~\ref{prop:four-quasigroups}).
Among these, laws that imply Sch-102 describe Boolean groups (rather than merely Abelian groups), as shown in Proposition~\ref{prop:102-boolean} below; their implication lattice is given in \autoref{fig:sublattice-102} on \autopageref{fig:sublattice-102}.
We distinguish between laws that do not require all bijections $\alpha,\beta,\gamma$ to be linear, presented in \S\ref{subsec:quasilinear}, those that require linearity, presented in \S\ref{subsec:noncomm}, and those that further require the bijections to commute in \S\ref{subsec:medial}, so that the whole structure describes a module over some ring.

\begin{remark}
  It is important to emphasize that the Abelian group structure (including the choice of zero) and the bijections are not necessarily uniquely defined, depending on the law.  For some of the equivalence classes below, an explicit construction of $(+,0,\alpha,\beta,\gamma)$ in terms of the quasigroup operations $\opA,\opB,\opC$ is given, but this has not be sought systematically.
\end{remark}

\begin{prop}[Sch-102 describes Boolean groups with a pair of bijections]
  \label{prop:102-boolean}
  Fix a quasigroup $(Q,\opA)$ obeying law Sch-102 $x \opA (y \opC z) = z \opA (y \opC x)$ and an element $0\in Q$.  Define
  \begin{equation}
    x + y \coloneqq ((0 \opC 0) \opB x) \opA (y \opC 0) , \qquad
    f(x) \coloneqq x \opA (0 \opC 0) , \qquad
    g(y) \coloneqq 0 \opA y .
  \end{equation}
  Then $f$ and $g$ are bijections, $(Q,+)$ is a Boolean group with neutral element~$0$, namely is Abelian and has $x+x=0$, and the quasigroup operations are expressed as
  \begin{equation}\label{xopy-102}
    x \opA y = f(x) + g(y) , \qquad
    x \opB y = f^{-1}(g(x)+y) , \qquad
    x \opC y = g^{-1}(x+f(y))
  \end{equation}
  Conversely, given $(Q,+,f,g)$ a Boolean group equipped with two bijections, the operations~\eqref{xopy-102} obey the quasigroup axioms~$\Qaxioms$ and law Sch-102.
\end{prop}

\begin{proof}
  Checking that~\eqref{xopy-102} obeys quasigroup axioms and Sch-102 is straightforward, so we concentrate on the other implication.
  Applying the four quasigroups theorem to Sch-102 gives that $x \opA y = \alpha(x) \times \mu(y)$ in some group $(Q,\times)$, not yet known to be Abelian.  Then by the quasigroup axiom $x \opA (y \opC x) = y$ we have $\mu(y \opC x) = (\alpha(x))^{-1} \times y$, where the $-1$ superscript denote the inverse in $(Q,\times)$.  Law Sch-102 then reads $\alpha(x) \times \mu(y \opC z) = \alpha(z) \times \mu(y \opC x)$ hence
  \begin{equation}
    \alpha(x) \times (\alpha(z))^{-1} \times y = \alpha(z) \times (\alpha(x))^{-1} \times y ,
  \end{equation}
  which simplifies to $a \times b^{-1} = b \times a^{-1}$ for all $a,b\in Q$ thanks to bijectivity of~$\alpha$.  Specializing $b$~to the neutral element~$e$ of $(Q,\times)$ we get $a=a^{-1}$, namely the group has exponent~$2$ and in particular is Abelian.  It is a Boolean group.

  We have $x \opA y = \alpha(x) \times \mu(y)$, $x \opB y = \alpha^{-1}(y \times \mu(x))$, and $x\opC y = \mu^{-1}(\alpha(y) \times x)$, where we omitted inversions as they are trivial in a Boolean group.
  There remains to evaluate, for some reference element~$0$ that is not necessarily the identity element~$e$ of~$(Q,\times)$, the operation~$+$:
  \begin{equation}
    \begin{aligned}
      x + y & = ((0 \opC 0) \opB x) \opA (y \opC 0) \\
      & = \alpha((0 \opC 0) \opB x) \times \mu(y \opC 0) \\
      & = x \times \mu(0 \opC 0) \times \alpha(0) \times y \\
      & = x \times y \times 0
    \end{aligned}
  \end{equation}
  as well as the bijections $f(x) = x \opA (0 \opC 0) = \alpha(x) \times \alpha(0) \times 0$ and $g(y) = 0 \opA y = \alpha(0) \times \mu(y)$ to check that
  \begin{equation}
    \begin{aligned}
    f(x) + g(y) & = \alpha(x) \times \alpha(0) \times 0 \times \alpha(0) \times \mu(y) \times 0 \\
      & = \alpha(x) \times \mu(y) = x \opA y .
    \end{aligned}
  \end{equation}
  This concludes the proof.
  As an alternative to this proof, the file \file{Sch-102.p} provides Prover9 code that successfully proves all the claims.
\end{proof}

The following two class lists concern general Abelian groups and Boolean groups, respectively.
When comparing laws that are parastrophes, one may need to change a bijection $f$ to~$f^{-1}$ and other similar redefinitions.
Given an Abelian group $(Q,+)$ and a bijection~$\alpha$, the operation $x\oplus y=\alpha^{-1}(\alpha(x)+\alpha(y))$ defines an Abelian group $(Q,\oplus)$; we use this freedom in writing Sch-54 in a more convenient manner than the expression $x \opA y = f(x)+Tf(y)$ that comes out of the four quasigroups theorem.
Some care is needed when considering conjunctions of equations, as the Abelian group operation~$+$ used in the description of each equation may not be the same.

\begin{classlist}\label{class:nonlinear-abelian}
Abelian groups twisted by bijections that are not necessarily linear.
\begin{itemize}

\item[\textbf{Sch-2}] $x \opA (y \opA z) \mkeqN{4362} y \opA (x \opA z)$ iff $Q$ is an \structure{Abelian group equipped with a bijection} $f:Q\to Q$ with $f(0)=0$; the operations are $x \opA y = f(x) + y$ and $x \opB y = f^{-1}(y - x)$ and $x \opC y = x - f(y)$.  The element $0$ is distinguished as the common value of all $x \opB x$, the addition as $x + y \coloneqq (0 \opB x) \opA y$, and $f(x) = x \opA 0$.

\item[\textbf{Sch-54.}] $x \opA (y \opB z) = y \opA (x \opB z)$ iff $Q$ is an \structure{affine $\ZZ[T]/(T^2+1)$-module equipped with a bijection $f:Q\to Q$}; the operations are $x \opA y = f^{-1}(x + Ty)$ and $x \opB y = -Tx + f(y)$ and $x \opC y = -Tf(x) + Ty$.  The operation $+$ and other data depend on a choice of origin $0\in Q$, from which one sets $f(x) = 0 \opB x$ and $Tx = 0 \opB (0 \opA x)$ and $x + y = 0 \opB (x \opA ((y \opA 0) \opC 0))$ and $-x = 0 \opB (0 \opA ((0 \opA 0) \opC x))$.  Equivalently, $x \mkeqN{562} y \opC (z \opC (y \opC (z \opC x)))$.

\item[Sch-100] $x \opA (y \opC z) = y \opA (x \opC z)$ iff $x \opA y = x-f(y)$ in an Abelian group equipped with a bijection~$f$ with $f(0)=0$ (para Sch-2).  Dubbed $D_2$ by Schr\"oder~\cite[p.641]{Schroder-1890}.

\item[Sch-140] $x \opB (y \opA z) = y \opB (x \opA z)$ iff $x \opA y = f(x-y)$ in an Abelian group equipped with a bijection~$f$ (para Sch-2).

\item[Sch-202] $x \opB (y \opC z) = y \opB (x \opC z)$ iff $x \opA y = -Tf(x) + Ty$ in an affine $\ZZ[T]/(T^2+1)$-module equipped with a bijection~$f$ (para Sch-54).

\item[Sch-224] $x \opC (y \opA z) = y \opC (x \opA z)$ iff $x \opA y = -Tx + f(y)$ in an affine $\ZZ[T]/(T^2+1)$-module equipped with a bijection~$f$ (para Sch-54).

\end{itemize}
\end{classlist}

\begin{classlist}\label{class:nonlinear-boolean}
Boolean groups twisted by bijections that are not necessarily linear.
\begin{itemize}
\item[\textbf{Sch-14}] $x \opA (y \opA z) = y \opA (z \opC x)$ iff $Q$ is a \structure{Boolean group equipped with a bijection} $f:Q\to Q$ with $f(0)=0$; the operation is $x \opA y = f(x) + y$.  This law is equivalent to Sch-2 together with the law $x+x=0$.  Equivalently, $x \opA y \mkeqN{3388} z \opA (x \opA (z \opA y))$.

\item[Sch-60] $x \opA (y \opB z) = y \opA (z \opC x)$ iff $x \opA y = f(x-y)$ in a Boolean group equipped with a bijection~$f$ with $f(0)=0$ (para Sch-14).

\item[Sch-101] $x \opA (y \opC z) = y \opA (z \opC x)$ iff $x \opA y = x-f(y)$ in a Boolean group equipped with a bijection~$f$ with $f(0)=0$ (para Sch-14).

\item[\textbf{Sch-102.}] $x \opA (y \opC z) = z \opA (y \opC x)$ iff $Q$ is a \structure{Boolean group equipped with two bijections} $f,g:Q\to Q$; the operation is $x \opA y = f(x)+g(y)$.  Cf.\ Proposition~\ref{prop:102-boolean}.

\item[\textbf{Sch-(54, 102)}] $x \opA (y \opB z) = y \opA (x \opB z)$ and $x \opA (y \opC z) = z \opA (y \opC x)$ iff $Q$ is an \structure{affine $(\ZZ/2\ZZ)[T]/(T^2+1)$-module equipped with a bijection $f:Q\to Q$}, equivalently a Boolean group equipped with a linear involution~$T$ and a bijection~$f$, with the operations $x \opA y = f^{-1}(x + Ty)$ and $x \opB y = Tx + f(y)$ and $x \opC y = Tf(x) + Ty$.  The data $(+,T,f)$ depend on a choice of origin $0\in Q$ as described for Sch-54 above.

\item[Sch-(102, 202)] $x \opA (y \opC z) = z \opA (y \opC x)$ and $x \opB (y \opC z) = y \opB (x \opC z)$ iff $x \opA y = Tf(x)+Ty$ in a Boolean group equipped with a linear involution~$T$ and a bijection~$f$ (para Sch-(54, 102)).

\item[Sch-(102, 224)] $x \opA (y \opC z) = z \opA (y \opC x)$ and $x \opC (y \opA z) = y \opC (x \opA z)$ iff $x \opA y = Tx + f(y)$ in a Boolean group equipped with a linear involution~$T$ and a bijection~$f$ (para Sch-(54, 102)).

\end{itemize}
\end{classlist}

\subsection{Linear noncommutative models (19 classes) (Class Lists~\ref{class:noncomm-abelian} and~\ref{class:noncomm-boolean})}
\label{subsec:noncomm}

We turn to classes whose operations read $x\opA y = Tx+Uy+c$ and $x\opB y = T^{-1}(-Ux+y-c)$ and $x\opC y = U^{-1}(x-Ty-c)$ in terms of an Abelian group $(Q,+)$, an element $c\in Q$, and two \emph{automorphisms} $T,U$ (bijective linear maps) that do not necessarily commute.  Classes for which $T,U$ commute are postponed to the next section.  As in the previous section we distinguish classes that only imply $(Q,+)$ is Abelian from those which restrict it to be a Boolean group.  This affects for instance the spectrum of finite cardinalities of quasigroups in these latter classes, as they must be powers of~$2$.

\begin{classlist}\label{class:noncomm-abelian}
Abelian groups twisted by linear bijections.  The operation $+$ and other data depend on a choice of origin $0\in Q$.
\begin{itemize}

\item[\textbf{Sch-(54, 202)}] $x \opA (y \opB z) = y \opA (x \opB z)$ and $x \opB (y \opC z) = y \opB (x \opC z)$ iff $Q$ is an \structure{Abelian group equipped with two linear anti-involutions~$T,U$ and an element~$c$}, namely $T^2=U^2=-1$.  The operations are $x \opA y = TUx+TUTy-TUc$ and $x \opB y = -Tx + UTy + c$ and $x \opC y = - TUTx + Ty - Tc$.  The linear anti-involutions (equivalently involutive anti-automorphisms) need not commute.

\item[Sch-(54, 224)] $x \opA (y \opB z) = y \opA (x \opB z)$ and $x \opC (y \opA z) = y \opC (x \opA z)$ (para Sch-(54, 202)).
\item[Sch-(202, 224)] $x \opB (y \opC z) = y \opB (x \opC z)$ and $x \opC (y \opA z) = y \opC (x \opA z)$ (para Sch-(54, 202)).

\item[\textbf{Sch-(54, 202, 224)}] $x \opA (y \opB z) = y \opA (x \opB z)$ and $x \opB (y \opC z) = y \opB (x \opC z)$ and $x \opC (y \opA z) = y \opC (x \opA z)$ iff $Q$ is an \structure{Abelian group equipped with two linear anti-commuting anti-involutions~$T,U$ and an element~$c$}, namely with $T^2=U^2=-1$ and $UT=-TU$, hence $TUT=U$ and $UTU=T$.  The operations are $x \opA y = TUx+Uy-TUc$ and $x \opB y = -Tx + UTy + c$ and $x \opC y = - Ux + Ty - Tc$.  The linear anti-involutions (equivalently involutive anti-automorphisms) need not commute.

\item[\textbf{Sch-(54, 408)}] $x \opA (y \opB z) = y \opA (x \opB z)$ and $x \opA (y \opC z) = (z \opC y) \opB x$ iff $Q$ is an \structure{Abelian group equipped with a linear anti-involution~$T$, a linear involution~$U$, and an element~$c$}, namely with $T^2=-1$ and $U^2=1$.  The operations are $x \opA y = Ux + UTy - Uc$ and $x \opB y = -Tx + Uy+c$ and $x \opC y = -TUx + Ty - Tc$.

\item[Sch-(54, 654)] $x \opA (y \opB z) = y \opA (x \opB z)$ and $x \opC (y \opB z) = (z \opB y) \opA x$ (para Sch-(54, 408)).
\item[Sch-(202, 408)] $x \opB (y \opC z) = y \opB (x \opC z)$ and $x \opA (y \opC z) = (z \opC y) \opB x$ (para Sch-(54, 408)).
\item[Sch-(202, 432)] $x \opB (y \opC z) = y \opB (x \opC z)$ and $x \opB (y \opA z) = (z \opA y) \opC x$ (para Sch-(54, 408)).
\item[Sch-(224, 432)] $x \opC (y \opA z) = y \opC (x \opA z)$ and $x \opB (y \opA z) = (z \opA y) \opC x$ (para Sch-(54, 408)).
\item[Sch-(224, 654)] $x \opC (y \opA z) = y \opC (x \opA z)$ and $x \opC (y \opB z) = (z \opB y) \opA x$ (para Sch-(54, 408)).

\item[Sch-(54, 408, 654)] $x \opA (y \opB z) = y \opA (x \opB z)$ and $x \opA (y \opC z) = (z \opC y) \opB x$ and $x \opC (y \opB z) = (z \opB y) \opA x$ iff $Q$ is an \structure{Abelian group equipped with linear maps $T,U$ with $T^2=-1$, $U^2=1$, $UT=-TU$, and an element~$c$}, with the operations $x \opA y = Ux + UTy - Uc$ and $x \opB y = -Tx + Uy+c$ and $x \opC y = -TUx + Ty - Tc$.  It is equivalent to Sch-(54, 408) with $UT=-TU$.

\item[Sch-(202, 408, 432)] $x \opB (y \opC z) = y \opB (x \opC z)$ and $x \opA (y \opC z) = (z \opC y) \opB x$ and $x \opB (y \opA z) = (z \opA y) \opC x$ (para Sch-(54, 408, 654)).

\item[Sch-(224, 432, 654)] $x \opC (y \opA z) = y \opC (x \opA z)$ and $x \opB (y \opA z) = (z \opA y) \opC x$ and $x \opC (y \opB z) = (z \opB y) \opA x$ (para Sch-(54, 408, 654)).
\end{itemize}
\end{classlist}

\begin{classlist}\label{class:noncomm-boolean}
Boolean groups twisted by linear bijections.
\begin{itemize}

\item[\textbf{Sch-(54, 102, 202)}] $x \opA (y \opB z) = y \opA (x \opB z)$ and $x \opA (y \opC z) = z \opA (y \opC x)$ and $x \opB (y \opC z) = y \opB (x \opC z)$ iff $Q$ is a \structure{Boolean group equipped with two linear involutions $T,U$ and an element~$c$} that is a fixed point of~$U$, with the operations $x \opA y = Ux + UTy + c$ and $x \opB y = Tx + Uy + c$ and $x \opC y = TUx + Ty + Tc$.  The involutions need not commute and $c$ need not vanish.  The data $(+,T,U,c)$ depend on a choice of origin $0\in Q$ as described for Sch-54 above.

\item[Sch-(54, 102, 224)] $x \opA (y \opB z) = y \opA (x \opB z)$ and $x \opA (y \opC z) = z \opA (y \opC x)$ and $x \opC (y \opA z) = y \opC (x \opA z)$ iff $x \opA y = TUx + Ty + Tc$ in a Boolean group equipped with two linear involutions $T,U$ and a fixed point $c$ of~$U$ (para Sch-(54, 102, 202)).

\item[Sch-(102, 202, 224)] $x \opA (y \opC z) = z \opA (y \opC x)$ and $x \opB (y \opC z) = y \opB (x \opC z)$ and $x \opC (y \opA z) = y \opC (x \opA z)$ iff $x \opA y = Tx + Uy + c$ in a Boolean group equipped with two linear involutions $T,U$ and a fixed point $c$ of~$U$ (para Sch-(54, 102, 202)).

\item[\textbf{Sch-(54, 264)}] $x \opA (y \opB z) = y \opA (x \opB z)$ and $x \opA (y \opA z) = (z \opA y) \opB x$ iff $Q$ is a \structure{Boolean group equipped with a linear $S_3$ action fixing an element~$c$}, or equivalently in terms of transpositions $T,U$ generating the permutation group~$S_3$, a Boolean group equipped with a pair of linear involutions $T,U$ with $(TU)^3=1$ and both fixing~$c$.  The operations are $x \opA y = Ux + UTy + c$ and $x \opB y = Tx + Uy + c$ and $x \opC y = TUx + Ty + c$.  The data $(+,T,U,c)$ depend on a choice of origin $0\in Q$ as described for Sch-54 above.

\item[Sch-(54, 504)] $x \opA (y \opB z) = y \opA (x \opB z)$ and $x \opB (y \opB z) = (z \opB y) \opC x$ iff $x \opA y = TUx + Ty + c$ in a Boolean group equipped with an $S_3$ action fixing~$c$ (para Sch-(54, 264)).

\item[Sch-(102, 336)] $x \opA (y \opC z) = z \opA (y \opC x)$ and $x \opA (y \opB z) = (z \opB y) \opB x$ iff $x \opA y = Tx + Uy + c$ in a Boolean group equipped with an $S_3$ action fixing~$c$ (para Sch-(54, 264)).
\end{itemize}
\end{classlist}

\subsection{Varieties of medial quasigroups (48 classes) (\autoref{class:medial})}
\label{subsec:medial}

% 3 4 7 8 9 10 15 20 25 34 38 56 79 2,18 2,54 2,56 4,29 4,56 4,100 4,101 4,102 4,183 4,224 4,408 25,54 29,56 34,54 34,56 34,93 34,102 34,140 34,432 54,100 54,127 56,102 56,202 56,654 103,202 127,140 161,224 4,29,56 4,102,224 34,54,102 54,202,408 54,224,654 56,102,202 202,224,432 54,102,202,224

To conclude our exploration of all \num{114} classes, we turn to classes whose models are commutative linear models on an Abelian group, or in other words affine modules over certain rings.
Specifically, we consider conjunctions of laws that imply the medial law
\begin{equation}
  (x \opA y) \opA (z \opA w) \mkeqN{934128} (x \opA z) \opA (y \opA w),
\end{equation}
equivalent to the same law for either of the other operations~$\opB$ or~$\opC$.  Then the Bruck-Murdoch-Toyoda theorem~\cite{Murdoch:1939,Toyoda:1941,Murdoch:1941,Bruck:1944} states that the quasigroup is isomorphic to some Abelian group $(Q,+)$ with a linear operation $x \opA y = A x + B y + c$ for some invertible linear maps $A,B\in\Aut(Q)$ that commute with each other, and an element $c\in Q$.

These data $(+,A,B,c)$ on~$Q$ are only defined up to a choice of reference element $0\in Q$.  (Nevertheless, the isomorphism class of the group $(Q,+)$ is uniquely determined, because isotopic groups are isomorphic.)
The other operations are also linear,
\[
x \opB y = - A^{-1}Bx + A^{-1}y - A^{-1}c , \qquad
x \opC y = B^{-1}x - B^{-1}Ay-B^{-1}c .
\]

\begin{remark}
If the linear map $A+B-1$ is invertible, the quasigroup has exactly one idempotent element: indeed, the linear equation $x = x\opA x = (A+B)x + c$ has exactly one solution.  In that case the abelian group structure $(Q,+)$ is made canonical by taking the reference element~$0$ to be this idempotent element.  In addition, $0 = 0 \opA 0 = (A+B)0 + c$ and linearity of $A,B$ implies that $c=0$.
\end{remark}

\begin{classlist}\label{class:medial}
\begin{itemize}

\item[\textbf{Sch-3}] $x \opA (y \opA z) = y \opA (z \opA x)$ iff $(Q,\opA)$ is an \structure{Abelian group}.
Dubbed $O_1$ by Schr\"oder~\cite[p.638]{Schroder-1890}.

\item[\textbf{Sch-4}] $x \opA (y \opA z) = z \opA (y \opA x)$ iff $Q$ is an \structure{affine $\ZZ[T,T^{-1}]$ module} with $x \opA y = T^2 x + Ty + c$.  The Abelian group is only defined up to a choice of reference element $0\in Q$: the subtraction $x - y \coloneqq (0 \opA x) \opC (y \opC (0 \opB 0))$ is such that $x-y=0$ and $x-(0-y)$ is an Abelian group operation.  Likewise, $T\colon x \mapsto (x \opC 0) \opA (0 \opA (0 \opB 0))$ depends on the choice of $0\in Q$.

\item[\textbf{Sch-7}] $x \opA (y \opA z) = y \opA (x \opB z)$ iff $Q$ is a \structure{Boolean group}, that is, all quasigroup operations are equal.  Equivalently, $x \mkeqN{895} y \opA ((x \opA z) \opA (y \opA z))$.

\item[Sch-8] $x \opA (y \opA z) = y \opA (z \opB x)$ iff $(Q,\opC)$ is an Abelian group (para Sch-3).

\item[\textbf{Sch-9}] $x \opA (y \opA z) = z \opA (x \opB y)$ iff $Q$ is an \structure{affine $\ZZ[T]/(T^2+1)$ module} with $x \opA y = -x+Ty+c$.  Equivalent to $x \mkeqN{546} y \opA (z \opA (x \opA (z \opA y)))$.
The module is defined up to a choice of $0\in Q$: the operation $x - y = (x\opA 0) \opA ((0 \opA 0) \opA y)$ gives the Abelian group subtraction, $B: x\mapsto (0 \opA 0) \opA x$, and $c = 0\opA 0$.

\item[\textbf{Sch-10}] $x \opA (y \opA z) = z \opA (y \opB x)$ iff $Q$ is an \structure{affine $(\ZZ/2\ZZ)[T]/(T^3-1)$ module} with $x \opA y = T x + T^2 y + c$, namely a Boolean group equipped with a group automorphism~$T$ of order~$3$ and some~$c$, defined up to a choice of origin $0\in Q$.  Equivalently, $x \opA y \mkeqN{3620} z \opA ((z \opA y) \opA x)$.  The quasigroup is semi-symmetric.

\item[Sch-15] $x \opA (y \opA z) = z \opA (x \opC y)$ iff $(Q,\opB)$ is an Abelian group (para Sch-3).

\item[Sch-20] $x \opA (y \opA z) = y \opB (z \opA x)$ iff $x \opA y = Tx - y + c$ with $T^2=-1$ (para Sch-9).

\item[Sch-25] $x \opA (y \opA z) = y \opB (x \opB z)$ iff $x \opA y = T (x + y) + c$ with $T^2=-1$ (para Sch-9).

\item[Sch-34] $x \opA (y \opA z) = z \opB (y \opC x)$ iff $x \opA y = Tx - T^{-1} y+c$ (para Sch-4).

\item[\textbf{Sch-38}] $x \opA (y \opA z) = y \opC (z \opA x)$ iff $Q$ is an \structure{affine $(\ZZ/2\ZZ)[T]/(T^7-1)$ module} with $x \opA y = Tx+T^3y+c$, namely a Boolean group equipped with a group automorphism~$T$ solving $(1+T)(1+T+T^3)(1+T^2+T^3)=0$, together with an element~$c$ with $(1+T)(1+T^2+T^3)c=0$.  The data $T,c$ and the operation~$+$ are determined up to a choice of $0\in Q$.  Equivalent to $x \mkeqN{2805} ((y \opA z) \opA (z \opA x)) \opA y$ (dual of \E{898}), or the same equation for $\opB$ or~$\opC$.

\item[Sch-56] $x \opA (y \opB z) = z \opA (y \opB x)$ iff $x \opA y = Tx + T^2 y + c$ (para Sch-4).

\item[Sch-79] $x \opA (y \opB z) = z \opB (x \opC y)$ iff $x \opA y = T^3x+Ty+c$ on a Boolean group with $T^7=1$ and $(1+T)(1+T^2+T^3)c=0$ (para Sch-38).

\item[\textbf{Sch-(2, 18)}] $x \opA (y \opA z) = y \opA (x \opA z)$ and $x \opA (y \opA z) = x \opB (z \opA y)$ iff $Q$ is a \structure{$(\ZZ/2\ZZ)[T]/T^2$ module} with $x \opA y = (1+T) x + y$, namely a Boolean group equipped with a nilpotent morphism $T: Q\to Q$ with $T^2=0$.

\item[\textbf{Sch-(2, 54)}] $x \opA (y \opA z) = y \opA (x \opA z)$ and $x \opA (y \opB z) = y \opA (x \opB z)$ iff $Q$ is a \structure{$\ZZ[T]/(T^2+1)$ module} with $x \opA y = T x + y$.  (In contrast to Sch-9 and its parastrophes this is a genuine module, not an affine one.)

\item[\textbf{Sch-(2, 56)}] $x \opA (y \opA z) = y \opA (x \opA z)$ and $x \opA (y \opB z) = z \opA (y \opB x)$ iff $Q$ is a \structure{$\ZZ[T]/(T^2-1)$ module} with $x \opA y = T x + y$.

\item[\textbf{Sch-(4, 29)}] $x \opA (y \opA z) = z \opA (y \opA x)$ and $x \opA (y \opA z) = x \opB (y \opC z)$ iff $Q$ is an \structure{affine $\ZZ[T]/(T^3-1)$ module} with $x \opA y = T^2 x - T y + c$.

\item[Sch-(4, 56)] $x \opA (y \opA z) = z \opA (y \opA x)$ and $x \opA (y \opB z) = z \opA (y \opB x)$ iff $x \opA y = T^2 x + T y + c$ with $T^3=1$ (para Sch-(4, 29)).

\item[Sch-(4, 100)] $x \opA (y \opA z) = z \opA (y \opA x)$ and $x \opA (y \opC z) = y \opA (x \opC z)$ iff $x \opA y = x + T y$ with $T^2=1$ (para Sch-(2, 56)).

\item[Sch-(4, 101)] $x \opA (y \opA z) = z \opA (y \opA x)$ and $x \opA (y \opC z) = y \opA (z \opC x)$ iff $x \opA y = x + (1+T) y$ on a Boolean group with $T$ a nilpotent morphism with $T^2=0$ (para Sch-(2, 18)).

\item[\textbf{Sch-(4, 102)}] $x \opA (y \opA z) = z \opA (y \opA x)$ and $x \opA (y \opC z) = z \opA (y \opC x)$ iff $Q$ is an \structure{affine $(\ZZ/2\ZZ)[T,T^{-1}]$ module} with $x \opA y = T^2 x + T y + c$.

\item[\textbf{Sch-(4, 183)}] $x \opA (y \opA z) = z \opA (y \opA x)$ and $x \opB (y \opB z) = x \opC (y \opA z)$ iff $Q$ is a \structure{$(\ZZ/2\ZZ)[T]/T^6$ module} with $x \opA y = (1+T^2) x + (1+T) y$.

\item[\textbf{Sch-(4, 224)}] $x \opA (y \opA z) = z \opA (y \opA x)$ and $x \opC (y \opA z) = y \opC (x \opA z)$ iff $Q$ is an \structure{affine $\ZZ[T]/(T^4+1)$ module} with $x \opA y = T^2 x + T y + c$.

\item[\textbf{Sch-(4, 408)}] $x \opA (y \opA z) = z \opA (y \opA x)$ and $x \opA (y \opC z) = (z \opC y) \opB x$ iff $Q$ is an \structure{affine $\ZZ[T]/(T^4-1)$ module} with $x \opA y = T^2 x + T y + c$.

\item[Sch-(25, 54)] $x \opA (y \opA z) = y \opB (x \opB z)$ and $x \opA (y \opB z) = y \opA (x \opB z)$ iff $x \opA y = (1+T)(x+y)$ on a Boolean group with $T^2=0$ (para Sch-(2, 18)).

\item[Sch-(29, 56)] $x \opA (y \opA z) = x \opB (y \opC z)$ and $x \opA (y \opB z) = z \opA (y \opB x)$ iff $x \opA y = (1+T) x + (1+T^2) y$ on a Boolean group with $T^6=0$ (para Sch-(4, 183)).

\item[Sch-(34, 54)] $x \opA (y \opA z) = z \opB (y \opC x)$ and $x \opA (y \opB z) = y \opA (x \opB z)$ iff $x \opA y = Tx - T^{-1} y+c$ with $T^4=-1$ (para Sch-(4, 224)).

\item[Sch-(34, 56)] $x \opA (y \opA z) = z \opB (y \opC x)$ and $x \opA (y \opB z) = z \opA (y \opB x)$ iff $x \opA y = - T^2 x + T y + c$ with $T^3=1$ (para Sch-(4, 29)).

\item[Sch-(34, 93)] $x \opA (y \opA z) = z \opB (y \opC x)$ and $x \opA (y \opB z) = x \opC (y \opC z)$ iff $x \opA y = (1+T)x + (1+T)^5 y$ on a Boolean group with $T^6=0$ hence $(1+T)^6=1$ (para Sch-(4, 183)).

\item[Sch-(34, 102)] $x \opA (y \opA z) = z \opB (y \opC x)$ and $x \opA (y \opC z) = z \opA (y \opC x)$ iff $x \opA y = T x + T^{-1} y + c$ with $T$ an automorphism of a Boolean group (para Sch-(4, 102)).

\item[Sch-(34, 140)] $x \opA (y \opA z) = z \opB (y \opC x)$ and $x \opB (y \opA z) = y \opB (x \opA z)$ iff $x \opA y = T(x-y)$ with $T^2=1$ (para Sch-(2, 56)).

\item[Sch-(34, 432)] $x \opA (y \opA z) = z \opB (y \opC x)$ and $x \opB (y \opA z) = (z \opA y) \opC x$ iff $x \opA y = T x - T^{-1} y + c$ with $T^4=1$ (para Sch-(4, 408)).

\item[Sch-(54, 100)] $x \opA (y \opB z) = y \opA (x \opB z)$ and $x \opA (y \opC z) = y \opA (x \opC z)$ iff $x \opA y = x + T y$ with $T^2=-1$ (para Sch-(2, 54)).

\item[\textbf{Sch-(54, 127)}] $x \opA (y \opB z) = y \opA (x \opB z)$ and $x \opA (y \opC z) = x \opC (y \opB z)$ iff $Q$ is a \structure{$\ZZ[T]/(T^4+1)$ module} with $x \opA y = Tx - T^3 y$.  (In contrast to Sch-(4, 224) and its parastrophes this is a genuine module, not an affine one.)

\item[Sch-(56, 102)] $x \opA (y \opB z) = z \opA (y \opB x)$ and $x \opA (y \opC z) = z \opA (y \opC x)$ iff $x \opA y = T x + T^2 y + c$ for some automorphism~$T$ of a Boolean group (para Sch-(4, 102)).

\item[Sch-(56, 202)] $x \opA (y \opB z) = z \opA (y \opB x)$ and $x \opB (y \opC z) = y \opB (x \opC z)$ iff $x \opA y = T x + T^2 y + c$ with $T^4=-1$ (para Sch-(4, 224)).

\item[Sch-(56, 654)] $x \opA (y \opB z) = z \opA (y \opB x)$ and $x \opC (y \opB z) = (z \opB y) \opA x$ iff $x \opA y = T x + T^2 y + c$ with $T^4=1$ (para Sch-(4, 408)).

\item[Sch-(103, 202)] $x \opA (y \opC z) = x \opB (y \opA z)$ and $x \opB (y \opC z) = y \opB (x \opC z)$ iff $x \opA y = Tx - T^2 y$ with $T^4=-1$ (para Sch-(54, 127)).

\item[Sch-(127, 140)] $x \opA (y \opC z) = x \opC (y \opB z)$ and $x \opB (y \opA z) = y \opB (x \opA z)$ iff $x \opA y = T(x-y)$ with $T^2=-1$ (para Sch-(2, 54)).

\item[Sch-(161, 224)] $x \opB (y \opA z) = x \opC (y \opB z)$ and $x \opC (y \opA z) = y \opC (x \opA z)$ iff $x \opA y = -T^2 x + T y$ with $T^4=-1$.

\item[\textbf{Sch-(4, 29, 56)}] $x \opA (y \opA z) = z \opA (y \opA x)$ and $x \opA (y \opA z) = x \opB (y \opC z)$ and $x \opA (y \opB z) = z \opA (y \opB x)$ iff $Q$ is an \structure{affine $(\ZZ/2\ZZ)[T]/(T^3-1)$ module} with $x \opA y = Tx + T^2x + c$.

\item[\textbf{Sch-(4, 102, 224)}] $x \opA (y \opA z) = z \opA (y \opA x)$ and $x \opA (y \opC z) = z \opA (y \opC x)$ and $x \opC (y \opA z) = y \opC (x \opA z)$ iff $Q$ is a \structure{$(\ZZ/2\ZZ)[T]/T^4$ module} with $x \opA y = (1+T^2)x + (1+T)y$

\item[Sch-(34, 54, 102)] $x \opA (y \opA z) = z \opB (y \opC x)$ and $x \opA (y \opB z) = y \opA (x \opB z)$ and $x \opA (y \opC z) = z \opA (y \opC x)$ iff $x \opA y = (1+T)(x+y)$ in a Boolean group equipped with a nilpotent endomorphism~$T$ with $T^4=0$ (para Sch-(4, 102, 224)).

\item[\textbf{Sch-(54, 202, 408)}] $x \opA (y \opB z) = y \opA (x \opB z)$ and $x \opB (y \opC z) = y \opB (x \opC z)$ and $x \opA (y \opC z) = (z \opC y) \opB x$ iff $Q$ is an \structure{affine $\ZZ[A,B]/(A^2-1,B^2+1)$ module} with $x \opA y = Ax+By+c$.

\item[Sch-(54, 224, 654)] $x \opA (y \opB z) = y \opA (x \opB z)$ and $x \opC (y \opA z) = y \opC (x \opA z)$ and $x \opC (y \opB z) = (z \opB y) \opA x$ iff $x \opA y = Ax+By+c$ with $A^2=-1$ and $B^2=1$ (para Sch-(54, 202, 408)).

\item[Sch-(56, 102, 202)] $x \opA (y \opB z) = z \opA (y \opB x)$ and $x \opA (y \opC z) = z \opA (y \opC x)$ and $x \opB (y \opC z) = y \opB (x \opC z)$ iff $x \opA y = (1+T)x + (1+T^2)y$ in a Boolean group with a nilpotent endomorphism~$T$ with $T^4=0$ (para Sch-(4, 102, 224)).

\item[Sch-(202, 224, 432)] $x \opB (y \opC z) = y \opB (x \opC z)$ and $x \opC (y \opA z) = y \opC (x \opA z)$ and $x \opB (y \opA z) = (z \opA y) \opC x$ iff $x \opA y = Ax+By+c$ with $A^2=B^2=-1$ (para Sch-(54, 202, 408)).

\item[\textbf{Sch-(54, 102, 202, 224)}] $x \opA (y \opB z) = y \opA (x \opB z)$ and $x \opA (y \opC z) = z \opA (y \opC x)$ and $x \opB (y \opC z) = y \opB (x \opC z)$ and $x \opC (y \opA z) = y \opC (x \opA z)$ iff $Q$ is a \structure{$(\ZZ/2\ZZ)[T,U]/(T^2,U^2)$ module} with $x \opA y = (1+T) x + (1+U)y$, namely a Boolean group equipped with two commuting nilpotent morphisms.

\end{itemize}
\end{classlist}

\section*{Acknowledgements}

I thank Stanley Burris for mentionning the work of Ernst Schr\"oder and giving me many details on the text, which is both in German and with long-unused terminology.  I also thank Jose Brox for related collaboration.

\bibliographystyle{plain}
\bibliography{Schroeder-refs}

\appendix

\setlength{\arraycolsep}{.5\arraycolsep}

\section{An interesting countermodel}
\label{app:counter}

The following seven-element quasigroup has all inverse properties, namely its sets of left multiplications and left divisions coincide, and likewise for middle and right, but the three corresponding involutions all differ pairwise:
\begin{equation}\label{quasigroup-all-IP}
\begin{array}{c|ccccccc}
  \copA & 0 & 1 & 2 & 3 & 4 & 5 & 6 \\
  \hline
  0 & 0 & 1 & 2 & 4 & 3 & 5 & 6 \\
  1 & 1 & 2 & 0 & 5 & 6 & 4 & 3 \\
  2 & 2 & 0 & 1 & 6 & 5 & 3 & 4 \\
  3 & 3 & 6 & 5 & 0 & 4 & 2 & 1 \\
  4 & 4 & 5 & 6 & 3 & 0 & 1 & 2 \\
  5 & 6 & 3 & 4 & 2 & 1 & 0 & 5 \\
  6 & 5 & 4 & 3 & 1 & 2 & 6 & 0
\end{array}
\quad
\begin{array}{c|ccccccc}
  \opB & 0 & 1 & 2 & 3 & 4 & 5 & 6 \\
  \hline
  0 & 0 & 1 & 2 & 3 & 4 & 6 & 5 \\
  1 & 2 & 0 & 1 & 5 & 6 & 4 & 3 \\
  2 & 1 & 2 & 0 & 6 & 5 & 3 & 4 \\
  3 & 3 & 6 & 5 & 4 & 0 & 1 & 2 \\
  4 & 4 & 5 & 6 & 0 & 3 & 2 & 1 \\
  5 & 5 & 4 & 3 & 2 & 1 & 0 & 6 \\
  6 & 6 & 3 & 4 & 1 & 2 & 5 & 0 \\
\end{array}
\quad
\begin{array}{c|ccccccc}
  \opC & 0 & 1 & 2 & 3 & 4 & 5 & 6 \\
  \hline
  0 & 0 & 2 & 1 & 3 & 4 & 5 & 6 \\
  1 & 1 & 0 & 2 & 6 & 5 & 4 & 3 \\
  2 & 2 & 1 & 0 & 5 & 6 & 3 & 4 \\
  3 & 4 & 6 & 5 & 0 & 3 & 1 & 2 \\
  4 & 3 & 5 & 6 & 4 & 0 & 2 & 1 \\
  5 & 5 & 3 & 4 & 2 & 1 & 6 & 0 \\
  6 & 6 & 4 & 3 & 1 & 2 & 0 & 5 \\
\end{array}
\end{equation}
The three involutions appearing in $f(x \opA y) = y \opA x$, $g(x \opB y) = y \opB x$, and $h(x \opC y) = y \opC x$ are (in cycle notation) $f=(34)(56)$ and $g=(12)(56)$ and $h=(12)(34)$.

\section{Pressure}
\label{sec:pressure}

Consider first the \textbf{empty class} $A=\emptyset$, whose variety consists of all quasigroups.
The total number of Cayley tables without imposing the invertibility conditions is~$n^{n^2}$.
The quasigroup condition amounts to focussing on Latin squares, only, whose number is known to be at least $(n!^2/n^n)^n = e^{n^2(\log n-2)+O(n\log n)}$.
This gives the desired asymptotics~\eqref{NnA-2var} with $p_\emptyset=0$.
Motivated by the fact that quasigroup conditions only reduce the number of valid Cayley tables by $e^{O(n^2)}$, a natural guess is that this also holds for other varieties.

Next, consider \textbf{commutativity} $x \opA y = y \opA x$, namely the symmetry of the Cayley table across the diagonal.  Only $n(n+1)/2$ entries are independent, which gives $n^{n(n+1)/2}$ entries.  Imposing the quasigroup conditions and quasigroup isomorphism should divide this number by $e^{O(n^2)}$ and $e^{O(n\log n)}$, respectively, which leads to my proposal that~\eqref{NnA-2var} would hold for this variety with $p_A=1/2$.

Concerning a \textbf{twist of commutativity} such as $x \opA y = f(y) \opA x$ for some bijection~$f$, note first that the number of commutative quasigroups gives a lower bound since they obey the equation with $f=\mathrm{id}$.  To place an upper bound, note that there are $n!$ possible permutations~$f$, and once this bijection is chosen the Cayley table can be filled one entry at a time.  Each time a non-diagonal entry $x \opA y$ is specified, the entry $f(y) \opA x$ must be filled too, and $f(x) \opA f(y)$, $f(f(y)) \opA f(x)$, and so all entries $f^k(x) \opA f^k(y)$ and $f^{k+1}(y) \opA f^k(x)$ must be filled, with $k\in\{0,\dots,K-1\}$, where $K$ is the greatest common divisor of the cycle lengths of~$f$ acting on $x$ and on~$y$.  This fixes two entries $x \opA y$ and $f(y) \opA x$ if $f(x)=x$ and $f(y)=y$, and more otherwise.  By filling the Cayley table in this way at most $n(n-1)/2$ off-diagonal entries must be specified, and $n$~diagonal ones.  This shows the desired upper bound, so \eqref{NnA-2var} would hold for $p_A=1/2$.  The remaining gap is that the lower bound has only been obtained by assuming that the quasigroup condition has an effect $e^{O(n^2)}$ only.

Next, consider \textbf{semi-symmetry} $x \opA y = x \opB y = x \opC y$.
Observe that from semi-symmetry and the relation between the operations, if $x \opA y = z$ then $z \opA x = (x \opB y) \opA x = y$ and $y \opA z = y \opA (x \opC y) = x$.
This easily provides an upper bound on the number of semi-symmetric quasigroups following the same strategy as for twisted commutativity above: fill in the Cayley table one entry at a time, with a choice $x \opA y = z$ leading to filling two other entries $z \opA x = y$ and $y \opA z = x$.  The accounting is made slightly trickier by diagonal entries, so the diagonal values of $x \opA x$ and $x \opB x$ and $x \opC x$ should be chosen first.  Altogether, there are $n^2/3+O(n)$ entries to be chosen at most.
A rigorous lower bound is provided by two distinct relations with Mendelsohn triple systems, namely decompositions of the complete directed graph on $n$ vertices into directed three-cycles.
\begin{itemize}
\item Semi-symmetric idempotent quasigroups~$Q$ are in one-to-one correspondence with Mendelsohn triple systems on~$Q$ by considering the directed three-cycles $(x,y,x\opA y)$ for $x,y$ distinct.  Indeed, idempotence $x \opA x = x$ implies that $x \opA y \neq x \opA x = x$ and $x \opA y\neq y\opA y = y$, so these are triples of distinct elements, and $y \opA (x \opA y) = x$ and $(x \opA y) \opA x = y$ so each directed pair appears in a single triple.
\item Semi-symmetric loops, namely semi-symmetric quasigroups~$Q$ with an identity~$e$, are in one-to-one correspondence with Mendelsohn triple systems on $Q\setminus\{e\}$ by considering the directed three-cycles $(x,y,x\opA y)$ for $x,y$ distinct and not equal to~$e$.
\end{itemize}
It is known~\cite{phelps1984number} that the number of $n$-element Mendelsohn triple systems is zero if $n\equiv 2\mod 3$ and otherwise equals $n^{n^2/3} e^{O(n^2)}$ at large $n\equiv 0,1\mod 3$.  The first identification shows the desired lower bound on $N_A(n)$ for $n\equiv 0,1\mod 3$ and the second identification for $n\equiv 1,2\mod 3$.
Altogether, for semi-symmetry, Conjecture~\ref{conj:NnA} is proven.

Twists of semi-symmetry work similarly to twists of commutativity.  A lower bound is given by semi-symmetric quasigroups themselves, and an upper bound by choosing first the needed bijection(s) and then filling in the Cayley table: each entry filled in will fix at least two other entries, or more depending on the cycle structure of the bijection.

Finally, \textbf{total symmetry} and its twists should be similar, with each $x \opA y = z$ with distinct $x,y,z$ fixing six entries $x \opA y$, $y \opA x$, $z \opA x$, $x \opA z$, $y \opA z$, $z \opA y$, which is why I expect $n^{n^2/6} e^{O(n^2)}$ quasigroups.  The lower bound can be partly proven by identifying totally symmetric idempotent quasigroups to Steiner systems, and totally symmetric loops to Steiner systems on $Q\setminus\{e\}$, but this only covers the cases $n\equiv 1,2,3,4\mod 6$ as Steiner systems only exist for $n\equiv 1,3\mod 6$.  Conjecture~\ref{conj:NnA} can nevertheless be proven for total symmetry\footnote{Total symmetry is \E{492} in the ETP numbering, and details on the spectrum of many such laws (for magmas, not only quasigroups) will appear as a separate publication.} for all congruences of~$n$.

An amusing question that comes out of this discussion would be to describe ``typical'' quasigroups in each class.  For instance, for a given implication $A\implies B$, what proportion of $B$-quasigroups (of size~$n$) are $A$-quasigroups.  I expect the proportion to tend to~$1$ as $n\to +\infty$ for many implications.

\section{Some interesting sublattices}
\label{sec:lattice-figs}

Here we collect a few interesting sublattices.

\begin{figure}[h]\centering
\begin{tikzpicture}[yscale=.7]
\node[draw](Sch-1) at (1.7883000000000002, -0.9) {\rotatebox{90}{$\copA$ commutes}};
\node[draw](Sch-88) at (1.0283, -4.5) {\rotatebox{90}{$\copA$ 2-TS}};
\node[draw](Sch-25) at (1.7883000000000002, -4.5) {\rotatebox{90}{Sch-25}};
\node[draw](Sch-60) at (2.5483000000000002, -4.5) {\rotatebox{90}{Sch-60}};
\node[draw](Sch-24) at (0.28832, -8.1) {\rotatebox{90}{TS}};
\node[draw](Sch-3) at (1.0283, -8.1) {\rotatebox{90}{$\copA$ Abel}};
\node[draw](Sch-25-54) at (1.8883, -8.1) {\rotatebox{90}{Sch-25-54}};
\node[draw](Sch-7) at (1.0283, -11.7) {\rotatebox{90}{Bool}};
\draw[-{stealth}](Sch-25.north)--(Sch-1.south);
\draw[-{stealth}](Sch-24.north)--(Sch-88.south);
\draw[-{stealth}](Sch-7.north)--(Sch-25-54.south);
\draw[-{stealth}](Sch-7.north)--(Sch-24.south);
\draw[-{stealth}](Sch-7.north)--(Sch-3.south);
\draw[-{stealth}](Sch-60.north)--(Sch-1.south);
\draw[-{stealth}](Sch-88.north)--(Sch-1.south);
\draw[-{stealth}](Sch-3.north)--(Sch-88.south);
\draw[-{stealth}](Sch-25-54.north)--(Sch-88.south);
\draw[-{stealth}](Sch-25-54.north)--(Sch-25.south);
\draw[-{stealth}](Sch-25-54.north)--(Sch-60.south);
\end{tikzpicture}
\caption{\label{fig:sublattice-1}Sublattice of classes that imply commutativity of~$\opA$ (law Sch-1).}
% python Schroeder_find_class.py hasse-tikz `python Schroeder_find_class.py stronger-short 1`
\end{figure}

\begin{figure}\centering
\begin{tikzpicture}[yscale=.9]
\node[draw](Sch-4) at (2.21, -0.9) {\rotatebox{90}{Sch-4}};
\node[draw](Sch-4-102) at (2.21, -4.5) {\rotatebox{90}{Sch-4-102}};
\node[draw](Sch-4-408) at (3.25, -4.5) {\rotatebox{90}{Sch-4-408}};
\node[draw](Sch-4-224) at (4.33, -4.5) {\rotatebox{90}{Sch-4-224}};
\node[draw](Sch-4-29) at (0.35997, -8.1) {\rotatebox{90}{Sch-4-29}};
\node[draw](Sch-4-56) at (1.26, -8.1) {\rotatebox{90}{Sch-4-56}};
\node[draw](Sch-4-183) at (2.21, -8.1) {\rotatebox{90}{Sch-4-183}};
\node[draw](Sch-4-100) at (3.2, -8.1) {\rotatebox{90}{$x \opC y = \sqrt{1} (x - y)$}};
\node[draw](Sch-4-102-224) at (4.38, -8.1) {\rotatebox{90}{Sch-4-102-224}};
\node[draw](Sch-9) at (5.43, -8.1) {\rotatebox{90}{Sch-9}};
\node[draw](Sch-4-29-56) at (1.26, -11.7) {\rotatebox{90}{Sch-4-29-56}};
\node[draw](Sch-15) at (2.23, -11.7) {\rotatebox{90}{$\opB$ Abel}};
\node[draw](Sch-3) at (2.97, -11.7) {\rotatebox{90}{$\copA$ Abel}};
\node[draw](Sch-4-101) at (3.83, -11.7) {\rotatebox{90}{Sch-4-101}};
\node[draw](Sch-10) at (1.55, -15.3) {\rotatebox{90}{Sch-10}};
\node[draw](Sch-7) at (2.51, -18.9) {\rotatebox{90}{Bool}};
\draw[-{stealth}](Sch-10.north)--(Sch-4-29-56.south);
\draw[-{stealth}](Sch-4-29-56.north)--(Sch-4-29.south);
\draw[-{stealth}](Sch-4-29-56.north)--(Sch-4-183.south);
\draw[-{stealth}](Sch-4-29-56.north)--(Sch-4-56.south);
\draw[-{stealth}](Sch-9.north)--(Sch-4-408.south);
\draw[-{stealth}](Sch-4-102-224.north)--(Sch-4-224.south);
\draw[-{stealth}](Sch-4-102-224.north)--(Sch-4-102.south);
\draw[-{stealth}](Sch-4-102-224.north)--(Sch-4-408.south);
\draw[-{stealth}](Sch-4-101.north)--(Sch-4-102-224.south);
\draw[-{stealth}](Sch-4-101.north)--(Sch-4-100.south);
\draw[-{stealth}](Sch-4-101.north)--(Sch-4-183.south);
\draw[-{stealth}](Sch-4-101.north)--(Sch-9.south);
\draw[-{stealth}](Sch-4-102.north)--(Sch-4.south);
\draw[-{stealth}](Sch-7.north)--(Sch-15.south);
\draw[-{stealth}](Sch-7.north)--(Sch-10.south);
\draw[-{stealth}](Sch-7.north)--(Sch-4-101.south);
\draw[-{stealth}](Sch-7.north)--(Sch-3.south);
\draw[-{stealth}](Sch-4-183.north)--(Sch-4-102.south);
\draw[-{stealth}](Sch-4-408.north)--(Sch-4.south);
\draw[-{stealth}](Sch-4-29.north)--(Sch-4.south);
\draw[-{stealth}](Sch-4-224.north)--(Sch-4.south);
\draw[-{stealth}](Sch-4-100.north)--(Sch-4-408.south);
\draw[-{stealth}](Sch-4-56.north)--(Sch-4.south);
\draw[-{stealth}](Sch-3.north)--(Sch-4-100.south);
\draw[-{stealth}](Sch-3.north)--(Sch-4-56.south);
\draw[-{stealth}](Sch-15.north)--(Sch-4-100.south);
\draw[-{stealth}](Sch-15.north)--(Sch-4-29.south);
\end{tikzpicture}
\caption{\label{fig:sublattice-4}Sublattice of classes that imply law Sch-4.}
% python Schroeder_find_class.py hasse-tikz `python Schroeder_find_class.py stronger-short 4`
\end{figure}

\begin{figure}\centering
\begin{tikzpicture}
\node[draw](Sch-102) at (9, -2) {\rotatebox{90}{Sch-102}};
\node[draw, text height=1.8cm](Sch-102-202) at (6.5, -4.5) {\rotatebox{90}{Sch-102-202}};
\node[draw, text height=1.8cm](Sch-54-102) at (7, -4.5) {\rotatebox{90}{Sch-54-102}};
\node[draw, text height=1.8cm](Sch-102-224) at (7.5, -4.5) {\rotatebox{90}{Sch-102-224}};
\node[draw, text height=1.8cm](Sch-56-102) at (9, -4.5) {\rotatebox{90}{Sch-56-102}};
\node[draw, text height=1.8cm](Sch-34-102) at (9.5, -4.5) {\rotatebox{90}{Sch-34-102}};
\node[draw, text height=1.8cm](Sch-4-102) at (10, -4.5) {\rotatebox{90}{Sch-4-102}};
\node[draw](Sch-38) at (12.4, -4.5) {\rotatebox{90}{Sch-38}};
\node[draw](Sch-79) at (13.2, -4) {\rotatebox{90}{Sch-79}};
\node[draw, text height=2.45cm](Sch-54-102-224) at (4.5, -8) {\rotatebox{90}{Sch-54-102-224}};
\node[draw, text height=2.45cm](Sch-102-202-224) at (5, -8) {\rotatebox{90}{Sch-102-202-224}};
\node[draw, text height=2.45cm](Sch-54-102-202) at (5.5, -8) {\rotatebox{90}{Sch-54-102-202}};
\node[draw, text height=1.2cm](Sch-14) at (7, -8) {\rotatebox{90}{Sch-14}};
\node[draw, text height=1.2cm](Sch-60) at (7.5, -8) {\rotatebox{90}{Sch-60}};
\node[draw, text height=1.2cm](Sch-101) at (8, -8) {\rotatebox{90}{Sch-101}};
\node[draw, text height=2.45cm](Sch-56-102-202) at (9, -8) {\rotatebox{90}{Sch-56-102-202}};
\node[draw, text height=2.45cm](Sch-34-54-102) at (9.5, -8) {\rotatebox{90}{Sch-34-54-102}};
\node[draw, text height=2.45cm](Sch-4-102-224) at (10, -8) {\rotatebox{90}{Sch-4-102-224}};
\node[draw](Sch-29-56) at (11, -8) {\rotatebox{90}{Sch-29-56}};
\node[draw](Sch-34-93) at (11.5, -8) {\rotatebox{90}{Sch-34-93}};
\node[draw](Sch-4-183) at (12, -8) {\rotatebox{90}{Sch-4-183}};
\node[draw](Sch-54-102-202-224) at (5, -11.5) {\rotatebox{90}{Sch-54-102-202-224}};
\node[draw, text height=1.8cm](Sch-54-264) at (5.8, -11) {\rotatebox{90}{Sch-54-264}};
\node[draw, text height=1.8cm](Sch-102-336) at (6.3, -11) {\rotatebox{90}{Sch-102-336}};
\node[draw, text height=1.8cm](Sch-54-504) at (6.8, -11) {\rotatebox{90}{Sch-54-504}};
\node[draw](Sch-4-29-56) at (12, -11) {\rotatebox{90}{Sch-4-29-56}};
\node[draw, text height=1.45cm](Sch-2-18) at (7.5, -15) {\rotatebox{90}{Sch-2-18}};
\node[draw, text height=1.45cm](Sch-25-54) at (8, -15) {\rotatebox{90}{Sch-25-54}};
\node[draw, text height=1.45cm](Sch-4-101) at (8.5, -15) {\rotatebox{90}{Sch-4-101}};
\node[draw](Sch-10) at (11, -14) {\rotatebox{90}{Sch-10}};
\node[draw](Sch-7) at (10.6, -17) {\rotatebox{90}{Bool}};
\draw[-{stealth}](Sch-54-264.north)--(Sch-54-102-224.south);
\draw[-{stealth}](Sch-2-18.north)--(Sch-56-102-202.south);
\draw[-{stealth}](Sch-2-18.north)--(Sch-54-264.south);
\draw[-{stealth}](Sch-2-18.north)--(Sch-54-102-202-224.south);
\draw[-{stealth}](Sch-2-18.north)--(Sch-14.south);
\draw[-{stealth}](Sch-2-18.north)--(Sch-29-56.south);
\draw[-{stealth}](Sch-54-102-224.north)--(Sch-102-224.south);
\draw[-{stealth}](Sch-54-102-224.north)--(Sch-54-102.south);
\draw[-{stealth}](Sch-56-102.north)--(Sch-102.south);
\draw[-{stealth}](Sch-38.north)--(Sch-102.south);
\draw[-{stealth}](Sch-54-102-202.north)--(Sch-102-202.south);
\draw[-{stealth}](Sch-54-102-202.north)--(Sch-54-102.south);
\draw[-{stealth}](Sch-10.north)--(Sch-4-29-56.south);
\draw[-{stealth}](Sch-102-202-224.north)--(Sch-102-224.south);
\draw[-{stealth}](Sch-102-202-224.north)--(Sch-102-202.south);
\draw[-{stealth}](Sch-54-102-202-224.north)--(Sch-102-202-224.south);
\draw[-{stealth}](Sch-54-102-202-224.north)--(Sch-54-102-202.south);
\draw[-{stealth}](Sch-54-102-202-224.north)--(Sch-54-102-224.south);
\draw[-{stealth}](Sch-101.north)--(Sch-102-224.south);
\draw[-{stealth}](Sch-4-102.north)--(Sch-102.south);
\draw[-{stealth}](Sch-54-102.north)--(Sch-102.south);
\draw[-{stealth}](Sch-102-224.north)--(Sch-102.south);
\draw[-{stealth}](Sch-60.north)--(Sch-54-102.south);
\draw[-{stealth}](Sch-4-101.north)--(Sch-4-102-224.south);
\draw[-{stealth}](Sch-4-101.north)--(Sch-4-183.south);
\draw[-{stealth}](Sch-4-101.north)--(Sch-101.south);
\draw[-{stealth}](Sch-4-101.north)--(Sch-54-102-202-224.south);
\draw[-{stealth}](Sch-4-101.north)--(Sch-54-504.south);
\draw[-{stealth}](Sch-29-56.north)--(Sch-56-102.south);
\draw[-{stealth}](Sch-7.north) to [bend right=20](Sch-38.south);
\draw[-{stealth}](Sch-7.north)--(Sch-10.south);
\draw[-{stealth}](Sch-7.north)--(Sch-25-54.south);
\draw[-{stealth}](Sch-7.north) to [bend right=20](Sch-79.south);
\draw[-{stealth}](Sch-7.north)--(Sch-2-18.south);
\draw[-{stealth}](Sch-7.north)--(Sch-4-101.south);
\draw[-{stealth}](Sch-14.north)--(Sch-102-202.south);
\draw[-{stealth}](Sch-54-504.north)--(Sch-54-102-202.south);
\draw[-{stealth}](Sch-34-102.north)--(Sch-102.south);
\draw[-{stealth}](Sch-102-336.north)--(Sch-102-202-224.south);
\draw[-{stealth}](Sch-34-54-102.north)--(Sch-34-102.south);
\draw[-{stealth}](Sch-34-54-102.north)--(Sch-54-102.south);
\draw[-{stealth}](Sch-4-29-56.north)--(Sch-4-183.south);
\draw[-{stealth}](Sch-4-29-56.north)--(Sch-34-93.south);
\draw[-{stealth}](Sch-4-29-56.north)--(Sch-29-56.south);
\draw[-{stealth}](Sch-79.north)--(Sch-102.south);
\draw[-{stealth}](Sch-4-102-224.north)--(Sch-102-224.south);
\draw[-{stealth}](Sch-4-102-224.north)--(Sch-4-102.south);
\draw[-{stealth}](Sch-4-183.north)--(Sch-4-102.south);
\draw[-{stealth}](Sch-102-202.north)--(Sch-102.south);
\draw[-{stealth}](Sch-34-93.north)--(Sch-34-102.south);
\draw[-{stealth}](Sch-56-102-202.north)--(Sch-102-202.south);
\draw[-{stealth}](Sch-56-102-202.north)--(Sch-56-102.south);
\draw[-{stealth}](Sch-25-54.north)--(Sch-102-336.south);
\draw[-{stealth}](Sch-25-54.north)--(Sch-60.south);
\draw[-{stealth}](Sch-25-54.north)--(Sch-34-54-102.south);
\draw[-{stealth}](Sch-25-54.north)--(Sch-54-102-202-224.south);
\draw[-{stealth}](Sch-25-54.north)--(Sch-34-93.south);
\end{tikzpicture}
\caption{\label{fig:sublattice-102}Sublattice of classes that imply law Sch-102 hence describe Boolean groups equipped with some structure (at most two bijections).  In particular the finite spectrum is a subset of the powers of two, and the number of magmas of size $n=2^k$ is at most~$n!^2$.}
% python Schroeder_find_class.py hasse-tikz `python Schroeder_find_class.py stronger-short 102`
\end{figure}

\begin{figure}\centering
\begin{tikzpicture}
\node[draw](Sch-) at (0,10) {\rotatebox{90}{no eq}};
\node[draw](Sch-3) at (-.6,-3) {\rotatebox{90}{$\copA$ Abel}};
\node[draw](Sch-15) at (0,-3) {\rotatebox{90}{$\opB$ Abel}};
\node[draw](Sch-8) at (.6,-3) {\rotatebox{90}{$\opC$ Abel}};
\node[draw](Sch-7) at (0,-5) {\rotatebox{90}{Bool}};
\node[draw](Sch-65) at (2.4,2) {\rotatebox{90}{$\copA$ Grp}};
\node[draw](Sch-105) at (3,2) {\rotatebox{90}{$\opB$ Grp}};
\node[draw](Sch-123) at (3.6,2) {\rotatebox{90}{$\opC$ Grp}};
\node[draw](Sch-432) at (-.6,8) {\rotatebox{90}{$\copA$ MIP}};
\node[draw](Sch-654) at (0,8) {\rotatebox{90}{$\opB$ MIP}};
\node[draw](Sch-408) at (.6,8) {\rotatebox{90}{$\opC$ MIP}};
\node[draw](Sch-34-140) at (-.6,1) {\rotatebox{90}{$x \copA y = \sqrt{1} (x - y)$}};
\node[draw](Sch-2-56) at (-0,1) {\rotatebox{90}{$x \opB y = \sqrt{1} (x - y)$}};
\node[draw](Sch-4-100) at (.6,1) {\rotatebox{90}{$x \opC y = \sqrt{1} (x - y)$}};
\node[draw](Sch-408-654) at (-.6,6) {\rotatebox{90}{$\copA$ IP}};
\node[draw](Sch-408-432) at (0,6) {\rotatebox{90}{$\opB$ IP}};
\node[draw](Sch-432-654) at (.6,6) {\rotatebox{90}{$\opC$ IP}};
\node[draw](Sch-408-432-654) at (0,4) {\rotatebox{90}{all IP}};
\draw[-{stealth}](Sch-432-654.north)--(Sch-654.south);
\draw[-{stealth}](Sch-432-654.north)--(Sch-432.south);
\draw[-{stealth}](Sch-105.north)--(Sch-408-432.south);
\draw[-{stealth}](Sch-408-654.north)--(Sch-654.south);
\draw[-{stealth}](Sch-408-654.north)--(Sch-408.south);
\draw[-{stealth}](Sch-408-432-654.north)--(Sch-432-654.south);
\draw[-{stealth}](Sch-408-432-654.north)--(Sch-408-654.south);
\draw[-{stealth}](Sch-408-432-654.north)--(Sch-408-432.south);
\draw[-{stealth}](Sch-2-56.north)--(Sch-408-432-654.south);
\draw[-{stealth}](Sch-7.north)--(Sch-8.south);
\draw[-{stealth}](Sch-7.north)--(Sch-3.south);
\draw[-{stealth}](Sch-7.north)--(Sch-15.south);
\draw[-{stealth}](Sch-654.north)--(Sch-.south);
\draw[-{stealth}](Sch-408-432.north)--(Sch-408.south);
\draw[-{stealth}](Sch-408-432.north)--(Sch-432.south);
\draw[-{stealth}](Sch-34-140.north)--(Sch-408-432-654.south);
\draw[-{stealth}](Sch-432.north)--(Sch-.south);
\draw[-{stealth}](Sch-123.north)--(Sch-432-654.south);
\draw[-{stealth}](Sch-15.north)--(Sch-4-100.south);
\draw[-{stealth}](Sch-15.north)--(Sch-105.south);
\draw[-{stealth}](Sch-15.north)--(Sch-34-140.south);
\draw[-{stealth}](Sch-4-100.north)--(Sch-408-432-654.south);
\draw[-{stealth}](Sch-8.north)--(Sch-34-140.south);
\draw[-{stealth}](Sch-8.north)--(Sch-2-56.south);
\draw[-{stealth}](Sch-8.north)--(Sch-123.south);
\draw[-{stealth}](Sch-408.north)--(Sch-.south);
\draw[-{stealth}](Sch-3.north)--(Sch-4-100.south);
\draw[-{stealth}](Sch-3.north)--(Sch-65.south);
\draw[-{stealth}](Sch-3.north)--(Sch-2-56.south);
\draw[-{stealth}](Sch-65.north)--(Sch-408-654.south);
\end{tikzpicture}
\caption{\label{fig:sublattice-associativity}Sublattice generated by associativity and its parastrophes.}
\end{figure}

\begin{figure}\centering
\begin{tikzpicture}
\node[draw](Sch-) at (3, 1) {\rotatebox{90}{no eq}};
\node[draw](Sch-432) at (2, -2) {\rotatebox{90}{$\copA$ MIP}};
\node[draw](Sch-654) at (2.6, -2) {\rotatebox{90}{$\opB$ MIP}};
\node[draw](Sch-408) at (3.2, -2) {\rotatebox{90}{$\opC$ MIP}};
\node[draw](Sch-56) at (5, -2) {\rotatebox{90}{Sch-56}};
\node[draw](Sch-4) at (6, -2) {\rotatebox{90}{Sch-4}};
\node[draw](Sch-127) at (0, -2) {\rotatebox{90}{$\copA$ MCIP}};
\node[draw](Sch-93) at (-1, -6) {\rotatebox{90}{$\copA$ MICIP}};
\node[draw](Sch-127-432) at (0, -6) {\rotatebox{90}{$\copA$ 4-TS}};
\node[draw](Sch-1) at (1, -6) {\rotatebox{90}{$\copA$ commutes}};
\node[draw](Sch-408-654) at (2, -6) {\rotatebox{90}{$\copA$ IP}};
\node[draw](Sch-408-432) at (2.6, -6) {\rotatebox{90}{$\opB$ IP}};
\node[draw](Sch-432-654) at (3.2, -6) {\rotatebox{90}{$\opC$ IP}};
\node[draw](Sch-2) at (4, -6) {\rotatebox{90}{Sch-2}};
\node[draw](Sch-336) at (1.5, -9) {\rotatebox{90}{$\copA$ EIP}};
\node[draw](Sch-408-432-654) at (2.6, -9) {\rotatebox{90}{all IP}};
\node[draw](Sch-4-56) at (4, -9) {\rotatebox{90}{Sch-4-56}};
\node[draw](Sch-56-654) at (6, -9) {\rotatebox{90}{Sch-56-654}};
\node[draw](Sch-100) at (7, -9) {\rotatebox{90}{Sch-100}};
\node[draw](Sch-4-408) at (8, -9) {\rotatebox{90}{Sch-4-408}};
\node[draw](Sch-88) at (1, -12) {\rotatebox{90}{$\copA$ 2-TS}};
\node[draw](Sch-65) at (2, -12) {\rotatebox{90}{$\copA$ Grp}};
\node[draw](Sch-2-56) at (6, -13) {\rotatebox{90}{$x \opB y = \sqrt{1} (x - y)$}};
\node[draw](Sch-4-100) at (7, -13) {\rotatebox{90}{$x \opC y = \sqrt{1} (x - y)$}};
\node[draw](Sch-3) at (3, -16) {\rotatebox{90}{$\copA$ Abel}};
\draw[-{stealth}](Sch-432-654.north)--(Sch-654.south);
\draw[-{stealth}](Sch-432-654.north)--(Sch-432.south);
\draw[-{stealth}](Sch-408-654.north)--(Sch-654.south);
\draw[-{stealth}](Sch-408-654.north)--(Sch-408.south);
\draw[-{stealth}](Sch-408-432-654.north)--(Sch-432-654.south);
\draw[-{stealth}](Sch-408-432-654.north)--(Sch-408-654.south);
\draw[-{stealth}](Sch-408-432-654.north)--(Sch-408-432.south);
\draw[-{stealth}](Sch-2-56.north)--(Sch-408-432-654.south);
\draw[-{stealth}](Sch-2-56.north)--(Sch-2.south);
\draw[-{stealth}](Sch-2-56.north)--(Sch-56-654.south);
\draw[-{stealth}](Sch-56.north)--(Sch-.south);
\draw[-{stealth}](Sch-93.north)--(Sch-127.south);
\draw[-{stealth}](Sch-432.north)--(Sch-.south);
\draw[-{stealth}](Sch-4-408.north)--(Sch-408.south);
\draw[-{stealth}](Sch-4-408.north)--(Sch-4.south);
\draw[-{stealth}](Sch-4-100.north)--(Sch-408-432-654.south);
\draw[-{stealth}](Sch-4-100.north)--(Sch-100.south);
\draw[-{stealth}](Sch-4-100.north)--(Sch-4-408.south);
\draw[-{stealth}](Sch-88.north)--(Sch-408-432-654.south);
\draw[-{stealth}](Sch-88.north)--(Sch-336.south);
\draw[-{stealth}](Sch-88.north)--(Sch-1.south);
\draw[-{stealth}](Sch-88.north)--(Sch-127-432.south);
\draw[-{stealth}](Sch-88.north)--(Sch-93.south);
\draw[-{stealth}](Sch-3.north)--(Sch-4-100.south);
\draw[-{stealth}](Sch-3.north)--(Sch-65.south);
\draw[-{stealth}](Sch-3.north)--(Sch-2-56.south);
\draw[-{stealth}](Sch-3.north)--(Sch-4-56.south);
\draw[-{stealth}](Sch-3.north)--(Sch-88.south);
\draw[-{stealth}](Sch-65.north)--(Sch-336.south);
\draw[-{stealth}](Sch-654.north)--(Sch-.south);
\draw[-{stealth}](Sch-127.north)--(Sch-.south);
\draw[-{stealth}](Sch-4-56.north)--(Sch-56.south);
\draw[-{stealth}](Sch-4-56.north)--(Sch-4.south);
\draw[-{stealth}](Sch-4-56.north)--(Sch-93.south);
\draw[-{stealth}](Sch-336.north)--(Sch-408-654.south);
\draw[-{stealth}](Sch-100.north)--(Sch-408.south);
\draw[-{stealth}](Sch-2.north)--(Sch-654.south);
\draw[-{stealth}](Sch-127-432.north)--(Sch-127.south);
\draw[-{stealth}](Sch-127-432.north)--(Sch-432.south);
\draw[-{stealth}](Sch-408-432.north)--(Sch-408.south);
\draw[-{stealth}](Sch-408-432.north)--(Sch-432.south);
\draw[-{stealth}](Sch-4.north)--(Sch-.south);
\draw[-{stealth}](Sch-1.north)--(Sch-432.south);
\draw[-{stealth}](Sch-56-654.north)--(Sch-56.south);
\draw[-{stealth}](Sch-56-654.north)--(Sch-654.south);
\draw[-{stealth}](Sch-408.north)--(Sch-.south);
\end{tikzpicture}
\caption{\label{fig:sublattice-group}Sublattice of consequences of the Abelian group law Sch-3, equivalent to commutativity and associativity of~$(Q,\opA)$.}
% python Schroeder_find_class.py hasse-tikz `python Schroeder_find_class.py weaker-short 3`
\end{figure}

\end{document}